\documentclass{article}

\usepackage{amsmath,amssymb}
\usepackage{amsthm}
\usepackage{amsbsy}
\usepackage{epsfig,graphics}
\usepackage{booktabs}
\usepackage{comment}
\usepackage{hyperref}
\usepackage{bbm}
\usepackage{cite}
\graphicspath{{./figures/}}

\usepackage{algorithm}
\usepackage{mathtools}
\usepackage{bm}
\usepackage{array,booktabs}
\usepackage{siunitx}

\usepackage{caption}

\renewcommand{\le}{\leqslant}

\newcommand{\R}{\ensuremath{\mathbb{R}}}

\newcommand{\Z}{\ensuremath{\mathbb{Z}}}

\newcommand{\F}{\ensuremath{\mathbf{F}}}
\newcommand{\f}{\ensuremath{\mathbf{f}}}
\newcommand{\x}{\ensuremath{\mathbf{x}}}
\newcommand{\uu}{\ensuremath{\mathbf{u}}}
\newcommand{\xii}{\ensuremath{\boldsymbol{\xi}}}
\newcommand{\X}{\ensuremath{\mathbf{X}}}

\newcommand{\dt}{\ensuremath{\Delta t}}

\newcommand\numberthis{\addtocounter{equation}{1}\tag{\theequation}}

\def\x{{\bf x}}

\def\dt{{\Delta t}}

\def\beq{\begin{equation}}
\def\eeq{\end{equation}}

\def\Z{{\bf Z}}

\def\rprime{\hbox{\hskip.25em\raise.5ex\hbox{$'$}\hskip.15em}}

\def\ddn1{{\frac{\partial}{\partial \nu_{\yb}}}}

\newcommand*{\affaddr}[1]{#1} 
\newcommand*{\affmark}[1][*]{\textsuperscript{#1}}

\def\nyu{Courant Institute of Mathematical Sciences, New York University,\\ New York, New York 10012}

\title{An Analysis of the Numerical Stability of the Immersed Boundary Method}

\author{%
Mengjian Hua\affmark[1], Charles S. Peskin\affmark[1]\\
$\,$ \\
\affaddr{\affmark[1]\nyu}\\
}
\date{\today}
\begin{document}

\maketitle

\begin{abstract} 
We present a numerical stability analysis of the immersed boundary(IB) method for a special case which is constructed so that Fourier analysis is applicable. We examine the stability of the immersed boundary method with the discrete Fourier transforms defined differently on the fluid grid and the boundary grid. This approach gives accurate theoretical results about the stability boundary since it takes the effects of the spreading kernel of the immersed boundary method on the numerical stability into account. In this paper, the spreading kernel is the standard 4-point IB delta function. A three-dimensional incompressible viscous flow and a no-slip planar boundary are considered. The case of a planar elastic membrane is also analyzed using the same analysis framework and it serves as an example of many possible generalizations of our theory. We present some numerical results and show that the observed stability behaviors are consistent with what are predicted by our theory.
\end{abstract}

\maketitle
\section{Introduction}

A large number of problems in biology are fluid-structure interaction problems and the immersed boundary method, originally introduced for the study of flow patterns around
 heart valves \cite{PESKIN1972252}, is both a mathematical formulation and a numerical method to treat such problems. 

A common difficulty encountered in the application of the immersed
 boundary method is numerical stiffness, requiring the use of small
 time steps, and this phenomenon has not been fully investigated.
 Some possible sources of the observed numerical stiffness are the
 singular nature of the force field applied by the immersed boundary
 to the fluid, and also the physical stiffness of the immersed
 boundary itself.
 
 Understanding the stability behavior of the immersed boundary method
 can provide guidance in adjusting numerical parameters
 and can help in the development of more stable IB schemes. The stability problem in which we are particularly interested concerns
 the use of target points in the modeling of no-slip boundaries. Here, immersed boundary points are held in place by stiff springs
 that connect them to target positions on the fixed no-slip boundary.
 This simple idea has been considerably generalized in the formulation
 of the penalty (pIB) immersed boundary method \cite{pIB}, which enables the
 simulation of immersed boundaries with mass, and the rigid pIB method \cite{pIB_rigid},
 which enables the simulation of immersed rigid bodies.  Although
 we do not study these generalizations of target points here, our results
 are probably applicable to them, at least qualitatively, since they
 involve the same kinds of spring-like forces employed in the same
 manner.

The model problem we choose is a special case in which Fourier analysis is applicable and we are able to make use of discrete Fourier transforms on the fluid grid and the boundary grid. To make Fourier analysis applicable, we assign some special features to the model problem we consider in this paper. These features are:
\begin{enumerate}
    \item[1)] Linearization of the Navier-Stokes equations by dropping the nonlinear terms
    \item[2)] Linearization of the boundary conditions by keeping the delta functions of the IB method centered at fixed locations.
    \item[3)] The use of uniform grids for discretization, with the boundary grid aligned parallel to the fluid grid, and with the meshwidths of the two grids related such that the fluid grid's meshwidth is an integer multiple of the boundary meshwidth. But note here that we do allow for an arbitrary translation of the boundary grid in relation to the fluid grid, so it is then possible that the two grids have no points in common. 
\end{enumerate}

Although these simplifications are needed to enable the stability analysis, we believe that the results are applicable, at least approximately, much more generally, and to some extent we have tested this and are reporting the numerical results in Section 3. 

It should be noted that linearization (2) of the boundary condition  is especially appropriate in the case of target points for modeling no-slip boundaries detailed in Section 2.2 and 2.3. In that case, we think that keeping delta functions fixed is something that should actually be done in practice. 

It should also be noted that linearizations (1) and (2) are appropriate for the study of small-amplitude vibrations of elastic membranes immersed in fluids, and that there is already a significant stiffness issue in this case, which is what we are addressing in Section 2.4 and Section 3.2. The linearization (1) and (2) can be formally justified by assuming that the fluid velocity is $O(\epsilon)$ and keeping only the lowest-order terms, although we do not give the details of that formal justification here. 

Moreover, the stability analysis can easily be generalized to the same type of immersed boundary that satisfies (3) but with different kinds of boundary forces. In other words, as long as the discretization of the boundary is the same, the form of the (spatially homogeneous) boundary force
 can be arbitrary, In this paper, we analyze the case of a planar elastic membrane as an example. 

Previous work closely related to the present paper has been done by
Stockie and Wetton \cite{stockie_1995} \cite{STOCKIE199941}.  They consider an elastic fiber immersed in
a two-dimensional, viscous, incompressible fluid, and they analyze the
small-amplitude modes of vibration of such a fiber to uncover the
source of the stiffness that is typically observed in immersed
boundary computations.  To make Fourier analysis applicable, they
consider the case in which the undisturbed configuration of the fiber
is straight.  A variety of time-stepping schemes are considered, and
the stability of each scheme is determined.  Spatial discretization is
not considered explicitly, but the effects of spatial discretization
are brought into the picture by restricting the modes under
consideration to those that can be represented on a grid of specified
meshwidth.  In the second paper cited above \cite{STOCKIE199941}, the smoothing
effect of the regularized delta function of the IB method is
considered, but in a continuous way.

The most important difference between our analysis and that of Stockie
and Wetton is that we do consider spatial discretization.  In fact,
the grids used to discretize the fluid and immersed boundary can be
different in our analysis, provided that the meshwidth used for the
fluid is an integer multiple of the boundary meshwidth, and the two
grids can be arbitrarily shifted with respect to each other, provided
that they are parallel.  Another difference is that our setting is
three-dimensional, and our domain is periodic.  Our focus is
specifically on the problem of a fixed boundary, modeled by immersed
boundary points held in place by stiff springs, but we also consider
an immersed membrane, and that problem is the direct generalization to
the 3D case of the 2D immersed fiber problem considered by Stockie and
Wetton.  Another difference is that we consider only one particular
time-stepping scheme, which can be described as a second-order
accurate Runge-Kutta scheme that is explicit in the immersed boundary
force (and in the nonlinear terms of the Navier-Stokes equations, but
these are not included in our analysis), but implicit in its handling
of the viscosity and incompressibility of the fluid.  This scheme
reduces to a leapfrog scheme for the particular problem that we
analyze, and the stability boundary that we find for this scheme turns
out to be independent of the fluid viscosity.

Stability analysis has also been done for a finite-element version
 of the immersed boundary method by Boffi et al. \cite{IBFE_boffi} and by Heltai \cite{HELTAI2008598}. The spatial discretization in these papers involves a variational formulation that avoids the explicit construction of a regularized
 delta function. The temporal discretization has in common with ours
 that the immersed boundary force is evaluated explicitly whereas the
 fluid solver is implicit, except for the nonlinear terms, which are
 not included in the analysis.

As additional background, we mention the proofs by Mori \cite{Mori_convergence} and by Liu
and Mori \cite{Liu2012PropertiesOD} \cite{Lpconvergence} of convergence of the IB method as applied to
problems such as the time-independent Stokes equations with a
force-field prescribed on an immersed boundary; and also the analysis
of the immersed boundary \textit{problem} by Mori et al. \cite{mori2019well} and by
Lin and Tong \cite{Lin_Tong}.  These works lay the foundation for a more complete
analysis of the immersed boundary method than anything attempted here.

The paper is organized as follows: in Section 2, we present our analytical derivations and results about the numerical stability of the immersed boundary method for the two model problems discussed above, i.e., the fixed boundary
and the elastic membrane, both immersed in a viscous incompressible fluid. Section 3 gives some numerical results that verify our stability analysis and provides numerical examples of the no-slip boundary and the elastic membrane. 

\section{Numerical stability analysis}

\subsection{Mathematical formulation}

We consider a cubic domain $ \mathbf{\Omega} = [0,L]\times[0,L]\times[0,L]$ with periodic boundary conditions and the domain of material coordinates $\mathbf{s} = (s_1,s_2)$ is $\mathbf{S} = [0,L]\times[0,L]$. The material coordinates $(s_1,s_2)$ are introduced to describe the position of the immersed boundary and its target position. Here, let $(x,y,z)$ be the Cartesian coordinates and let the boundary be a plane parallel to the xy-plane and located at $z = \sigma_3 \in [0,L]$. Ideally, we want the immersed boundary to be no-slip and therefore the immersed boundary is expected to be fixed at its initial position, which we call the target position here. Thus, the target position, denoted here by $\X^0$, in the Cartesian coordinates is given by $\X^0(s_1,s_2) = \X(s_1,s_2,0) = (s_1 + \sigma_1,s_2 + \sigma_2,\sigma_3)$, where $\X(s_1,s_2,t)$ is the boundary position and $t$ is the time. Note that the immersed boundary $\X(s_1,s_2,t)$ is assumed to move with the fluid near it by the formulation of the immersed boundary method and we cannot strictly fix it at its target position $\X^0(s_1,s_2)$. Therefore, the immersed boundary $\X(s_1,s_2,t)$ may change its position over time and it depends on time $t$. Equation \eqref{feedback_force}, which we will introduce later in this section, describes how we exert a spring-like feedback force on the immersed boundary $\X(s_1,s_2,t)$ in order to keep it near its initial positions $\X^0(s_1,s_2)$. 

The equations of motion are
\begin{equation} \label{1}
    \rho \frac{\partial \mathbf{u}}{\partial t} + \nabla \mathbf{p} = \mu \Delta \mathbf{u} + \mathbf{f}
\end{equation}
\begin{equation} \label{2}
    \nabla \cdot \mathbf{u} = 0
\end{equation}
\begin{equation}\label{3}
    \mathbf{f}(\mathbf{x},t) = \int_\mathbf{S} \mathbf{F}(s_1,s_2,t)\delta(\mathbf{x} - \mathbf{X}^0 (s_1,s_2)) d s_1 d s_2
\end{equation}
\begin{equation}\label{4}
    \frac{\partial \mathbf{X}}{\partial t} (s_1,s_2,t)= \int_{\R^3} \mathbf{u} (\mathbf{x},t)\delta(\mathbf{x} - \mathbf{X}^0 (s_1,s_2)) d \mathbf{x}
\end{equation}
\begin{equation}\label{feedback_force}
    \mathbf{F} (s_1,s_2,t) = -K (\mathbf{X}(s_1,s_2,t) - \mathbf{X}^0(s_1,s_2))
\end{equation}

Equations \eqref{1} and \eqref{2} are the time-dependent incompressible Stokes equations, which are the linearized incompressible Navier-Stokes equations with the nonlinear convective term being dropped. Equations \eqref{3} and \eqref{4} are interaction equations which translate between Eulerian and Lagrangian variables. Equation \eqref{3} describes the body force $\mathbf{f}(\mathbf{x},t)$ that is applied to the
fluid and equation \eqref{4} describes the evaluation of the fluid velocity at the immersed boundary. As mentioned above, equation \eqref{feedback_force} describes the force generated by the stiff springs that hold the immersed boundary $\mathbf{X}(s_1,s_2,t)$ in place.

Since the velocity evaluated at the boundary in \eqref{4} should be zero because of the no-slip boundary condition, $\mathbf{X}(s_1,s_2,t) - \mathbf{X}^0(s_1,s_2)$ is the integral over time of the error that has been made in enforcing the no-slip condition.  That is why our formula for $\mathbf{F}(s_1,s_2,t)$ in equation \eqref{feedback_force} can be interpreted as a feedback mechanism (specifically, an integral controller) for enforcing the no-slip condition.

\subsection{Numerical Scheme}
The fluid domain is discretized by an $N\times N \times N$ grid with a uniform mesh width $h = L/N$. Let $Z_N = \{0,1,\cdots,N-1\}$ and $Z_{NP} = \{0,1,\cdots,NP-1\}$ be sets of integer indices, where $P\in \mathbb{N}$. Then, the fluid grid is defined by $ \mathbf{\Omega}_h = \{\mathbf{j}h: \mathbf{j} = (j_1,j_2,j_3) \in Z_N^3 \}$ and the grid points are $\mathbf{x}_{\mathbf{j}} = (j_1 h,j_2 h,j_3 h)$. Let the boundary grid be defined by $ \mathbf{S}_{h_B} = \{\mathbf{k} h_B + \bm{\sigma} : \mathbf{k} = (k_1,k_2) \in Z_{NP}^2 \}$, where $h_B = L/(NP)$ is the mesh width of the boundary grid and $\boldsymbol{\sigma} = (\mathbf{\sigma}_1,\mathbf{\sigma}_2,\mathbf{\sigma}_3)$ is a three-dimensional shift of the immersed boundary. Note that the boundary grid plane is defined such that it is parallel to the xy-plane and its z-coordinate, namely $\sigma_3$, is general. In other words, the boundary plane does not have to lie on one of the fluid grid planes. The target position of the discretized immersed boundary is thus given by $\mathbf{X}^0_{k_1,k_2} = \mathbf{X}^0 (k_1 h_B,k_2 h_B) = (k_1 h_B + \mathbf{\sigma}_1,k_2 h_B + \mathbf{\sigma}_2,\mathbf{\sigma}_3)$ and the discretized boundary itself is given by $\X_{k_1,k_2} (t) = \X(k_1 h_B,k_2 h_B, t)$. Similarly, the force at the boundary $\F_{k_1,k_2} (t)$ is defined by 
$\F_{k_1,k_2} (t) = \F(k_1 h_B,k_2 h_B, t)$. In the above definitions, we use subscripts to indicate indices $(k_1,k_2) \in Z_{NP}^2$.  
The spatial discretization of equation (1) - (5) is then as follows
\begin{equation} 
    \rho \frac{\partial \mathbf{u}}{\partial t} + \nabla_h \mathbf{p} = \mu \Delta_h \mathbf{u} + \mathbf{f}
\end{equation}
\begin{equation} 
    \nabla_h \cdot \mathbf{u} = 0
\end{equation}
\begin{equation} \label{8}
    \mathbf{f}(\mathbf{x}_{\mathbf{j}},t) = \sum_{\mathbf{k}\in Z_{NP}^2} \mathbf{F}_{k_1,k_2}(t) \delta_h (\mathbf{x}_{\mathbf{j}} - \mathbf{X}^0_{k_1,k_2}) h_B^2
\end{equation}
\begin{equation}
    \frac{\partial \mathbf{X}_{k_1,k_2}}{\partial t} (t)= \sum_{\mathbf{j}\in Z_N^3} \mathbf{u}( \mathbf{x}_{\mathbf{j}},t)\delta_h (\mathbf{x}_{\mathbf{j}} - \mathbf{X}^0_{k_1,k_2}) h^3
\end{equation}
\begin{equation}\label{10}
    \mathbf{F}_{k_1,k_2} (t) = -K (\mathbf{X}_{k_1,k_2}(t) - \mathbf{X}^0_{k_1,k_2})
\end{equation}

Note that the smoothed Dirac delta function $\delta_h$ is evaluated at target positions $\mathbf{X}^0_{k_1,k_2}$ instead of $ \mathbf{X}_{k_1,k_2} (t)$. For a sufficiently large $K$, we expect $\mathbf{X}_{k_1,k_2}(t)$ to be approximately equal to $\mathbf{X}^0_{k_1,k_2}$, and by using $\mathbf{X}^0_{k_1,k_2}$, we linearizes the boundary condition. 
The difference operators $\nabla_h$ and $\Delta_h$ are defined as follows:
\begin{equation}
    (\nabla_h p)_\alpha (\mathbf{x}_{\mathbf{j}}) = \frac{p(\mathbf{x}_{\mathbf{j}} + h \mathbf{e}_\alpha) - p(\mathbf{x}_{\mathbf{j}} - h\mathbf{e}_\alpha)}{2h},\quad \alpha = 1,2,3
\end{equation}

\begin{equation}
     (\nabla_h \cdot \mathbf{u})(\mathbf{x}_{\mathbf{j}}) = \sum_{\alpha =1}^3 \frac{\mathbf{u}(\mathbf{x}_{\mathbf{j}} + h \mathbf{e}_\alpha) - \mathbf{u}(\mathbf{x}_{\mathbf{j}} - h\mathbf{e}_\alpha)}{2h}
\end{equation}
\begin{equation}
     (\Delta_h \mathbf{u})(\mathbf{x}_{\mathbf{j}}) = \sum_{\alpha =1}^3 \frac{\mathbf{u}(\mathbf{x}_{\mathbf{j}} + h \mathbf{e}_\alpha) + \mathbf{u}(\mathbf{x}_{\mathbf{j}} - h\mathbf{e}_\alpha) - 2\mathbf{u}(\mathbf{x}_{\mathbf{j}})}{h^2}
\end{equation}

In the above definitions, $\{\mathbf{e}_1,\mathbf{e}_2,\mathbf{e}_3\}$ is the standard basis of $\R^3$ and $\mathbf{j} \in Z_N^3$, so $\x_{\mathbf{j}} \in \Omega_h$.

We can eliminate $\mathbf{X}_{k_1,k_2}(t)$ from our system by differentiating \eqref{10} with respect to $t$ and it follows that

\begin{equation}\label{14}
    \frac{\partial \mathbf{F}_{k_1,k_2}}{\partial t} (t)= -K \sum_{\mathbf{j}\in Z_N^3} \mathbf{u} (\mathbf{x}_\mathbf{j},t)\delta(\mathbf{x}_\mathbf{j} - \mathbf{X}^0_{k_1,k_2}) h^3
\end{equation}

Note that $\F(k_1,k_2,t)$ is then proportional to the accumulated
 error up to time t that has occurred in enforcing the no-slip
 condition at the location $\X^0(k_1,k_2)$. In the above equations, the smoothed Dirac delta function $\delta_h$ is the standard IB 4-point delta function defined \cite{peskin_2002} by 
\begin{equation}
    \delta_h(\mathbf{x}_\mathbf{j}) = \frac{1}{h^3} \phi(j_1)\phi(j_2)\phi(j_3)
\end{equation}
where 
\begin{equation}
    \phi(r) = \begin{cases}
    \frac{3-2|r|+ \sqrt{1+4|r|-4 r^2}}{8}, \quad |r| \leq 1\\
    \frac{5-2|r|- \sqrt{-7+12|r|-4 r^2}}{8}, \quad |r| \in [1,2]\\
    0, \quad |r|>2
    \end{cases}
\end{equation}
\begin{figure}[h]
\centering
\includegraphics[width=0.9\textwidth]{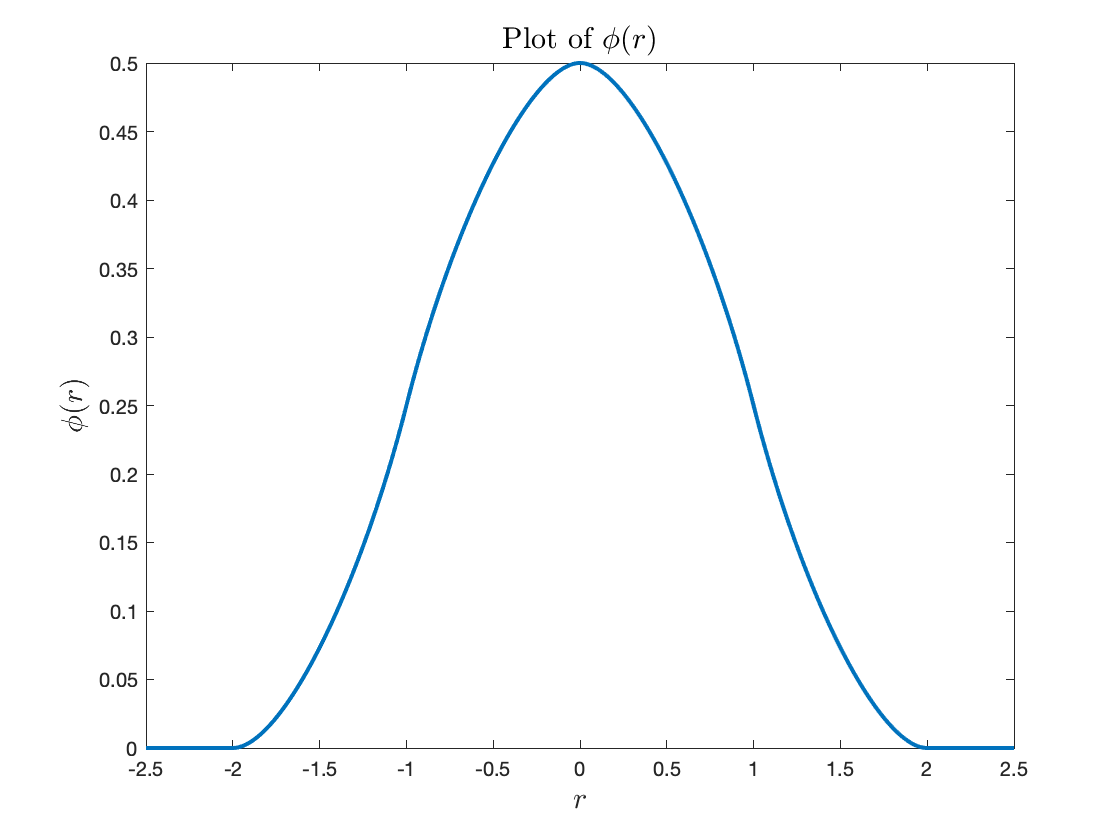}
\caption{A plot of $\phi(r)$, which is a bell-shaped kernel used to spread force fields and interpolate velocity fields in the immersed boundary method}
\label{fig:plot_IB4}
\end{figure}

Let $\boldsymbol\epsilon = \boldsymbol\sigma/h$, so that $\boldsymbol\epsilon$ is the shift vector in
units of meshwidth of the fluid grid. Then, equation \eqref{8} and \eqref{14} can be rewritten as follows
\begin{equation} 
    \mathbf{f}(\mathbf{x}_\mathbf{j},t) = \sum_{\mathbf{k}\in Z_{NP}^2} \mathbf{F}_{k_1,k_2}(t) \phi(j_1 - \frac{k_1}{P}-\epsilon_1) \phi(j_2 - \frac{k_2}{P} -\epsilon_2) \phi(j_3  -\epsilon_3)\frac{h_B^2}{h^3}
\end{equation}
\begin{equation}
    \frac{\partial \mathbf{F}_{k_1,k_2}}{\partial t} (t)= -K \sum_{\mathbf{j}\in Z_N^3} \mathbf{u} (\mathbf{x}_\mathbf{j},t) \phi(j_1 - \frac{k_1}{P}-\epsilon_1) \phi(j_2 - \frac{k_2}{P} -\epsilon_2) \phi(j_3  -\epsilon_3)
\end{equation}
For discretization in time, we first write down a second-order Runge-Kutta scheme as in \cite{peskin_2002}, but this scheme will turn out to be equivalent to a simple leapfrog scheme because we have here dropped the nonlinear convective term of the Navier-Stokes equations and also because the centering of the smoothed Dirac delta function $\delta_h$ at fixed positions. 

In writing this scheme, it is helpful to introduce the following notation
\begin{equation}\label{big_U}
    \mathbf{U}_{k_1,k_2} (t) = \sum_{\mathbf{j}\in Z_N^3} \mathbf{u} (\mathbf{x}_\mathbf{j},t) \phi(j_1 - \frac{k_1}{P}-\epsilon_1) \phi(j_2 - \frac{k_2}{P} -\epsilon_2) \phi(j_3  -\epsilon_3)
\end{equation}

From now on, we use superscripts to indicate the time step. The Runge-Kutta scheme is then as follows
\begin{equation}\label{force_first}
    \mathbf{F}^{n+\frac{1}{2}}_{k_1,k_2} = \F^n_{k_1,k_2} - \frac{\Delta t}{2} K \mathbf{U}^{n}_{k_1,k_2}
\end{equation}
\begin{equation}
    \f^{n+\frac{1}{2}} (\x_\mathbf{j}) = \sum_{\mathbf{k}\in Z_{NP}^2} \mathbf{F}^{n+\frac{1}{2}}_{k_1,k_2} \phi(j_1 - \frac{k_1}{P}-\epsilon_1) \phi(j_2 - \frac{k_2}{P} -\epsilon_2) \phi(j_3  -\epsilon_3)\frac{h_B^2}{h^3}
\end{equation}
\begin{equation}\label{sys1}
\begin{dcases*}
\rho \frac{\uu^{n+\frac{1}{2}} - \uu^n }{\dt/2} + \nabla_h \tilde{\mathbf{p}}^{n+\frac{1}{2}} = \mu \Delta_h \mathbf{u}^{n+\frac{1}{2}} + \mathbf{f}^{n+\frac{1}{2}} \\
        \nabla_h \cdot \mathbf{u}^{n+\frac{1}{2}} = 0
\end{dcases*}
\end{equation}
\begin{equation}\label{sys2}
\begin{dcases}
\rho \frac{\uu^{n+1} - \uu^n }{\dt} + \nabla_h \mathbf{p}^{n+\frac{1}{2}} = \mu \Delta_h \frac{\uu^{n+1} + \uu^n }{2} + \mathbf{f}^{n+\frac{1}{2}}\\
    \nabla_h  \cdot \mathbf{u}^{n+1} = 0
\end{dcases}
\end{equation}
\begin{equation}\label{force_last}
    \mathbf{F}^{n+1}_{k_1,k_2} = \F^n_{k_1,k_2} - \dt K \mathbf{U}^{n+\frac{1}{2}}_{k_1,k_2}
\end{equation}

To see how this scheme simplifies, we notice that systems of equations \eqref{sys1} and \eqref{sys2}, which have the same source term, are independent of each other, since $\uu^{n+\frac{1}{2}}$ does not appear in \eqref{sys2} as it would if we were considering the Naiver-Stokes equations. Moreover, if we have the solution of \eqref{sys2}, then we also have the solution of \eqref{sys1} simply by setting 
\begin{equation}\label{simp_leap}
\begin{dcases}
    \uu^{n+\frac{1}{2}} = \frac{\uu^n + \uu^{n+1}}{2}\\
   \tilde{\mathbf{p}}^{n+\frac{1}{2}} = \mathbf{p}^{n+\frac{1}{2}}
\end{dcases}
\end{equation}
Equation \eqref{simp_leap} implies that
\begin{equation}
    \mathbf{U}^{n+\frac{1}{2}} = \frac{\mathbf{U}^n + \mathbf{U}^{n+1}}{2}
\end{equation}
But note that it does so only because of the use of delta functions
 centered at fixed locations.
Thus, we can dispense with \eqref{sys1}, and also we can rewrite \eqref{force_last} as 
\begin{equation}
    \mathbf{F}^{n+1}_{k_1,k_2} = \mathbf{F}^{n+\frac{1}{2}}_{k_1,k_2} - \frac{\dt}{2} K \mathbf{U}^{n+1}_{k_1,k_2}
\end{equation}
Lowering n by 1 gives 
\begin{equation}
    \mathbf{F}^{n}_{k_1,k_2} = \mathbf{F}^{n-\frac{1}{2}}_{k_1,k_2} - \frac{\dt}{2} K \mathbf{U}^{n}_{k_1,k_2}
\end{equation}
and substitution of this into \eqref{force_first} gives 
\begin{equation}
    \mathbf{F}^{n+\frac{1}{2}}_{k_1,k_2} = \mathbf{F}^{n-\frac{1}{2}}_{k_1,k_2} - \dt K \mathbf{U}^{n}_{k_1,k_2}
\end{equation}
Therefore, we are left with a leapfrog scheme that can be written as follows:
\begin{equation}\label{leapfrog_1}
    \mathbf{F}^{n+\frac{1}{2}}_{k_1,k_2} = \mathbf{F}^{n-\frac{1}{2}}_{k_1,k_2} - \dt K \mathbf{U}^{n}_{k_1,k_2}
\end{equation}
\begin{equation}\label{leapfrog_2}
    \f^{n+\frac{1}{2}} (\x_\mathbf{j}) = \sum_{\mathbf{k}\in Z_{NP}^2} \mathbf{F}^{n+\frac{1}{2}}_{k_1,k_2} \phi(j_1 - \frac{k_1}{P}-\epsilon_1) \phi(j_2 - \frac{k_2}{P} -\epsilon_2) \phi(j_3  -\epsilon_3)\frac{h_B^2}{h^3}
\end{equation}
\begin{equation}\label{leapfrog_3}
\begin{dcases}
\rho \frac{\uu^{n+1} - \uu^n }{\dt} + \nabla_h \mathbf{p}^{n+\frac{1}{2}} = \mu \Delta_h \frac{\uu^{n+1} + \uu^n }{2} + \mathbf{f}^{n+\frac{1}{2}}\\
    \nabla_h \cdot \mathbf{u}^{n+1} = 0
\end{dcases}
\end{equation}
\subsection{Stability analysis of the scheme}
To analyze the numerical stability of the above leapfrog scheme, we introduce the discrete Fourier transform defined on the fluid grid by 
\begin{equation}\label{fourier_transform_1}
    \widehat{\uu}(\mathbf{\xii}) = \sum_{\mathbf{j} \in Z_N^3} e^{-i \frac{2\pi}{N} \mathbf{j} \cdot \mathbf{\xii}} \uu(\mathbf{x}_\mathbf{j})
\end{equation}
\begin{equation}\label{inv_Fourier_1}
    \uu(\mathbf{x}_\mathbf{j}) = \frac{1}{N^3} \sum_{\mathbf{\xii} \in Z_N^3} e^{i \frac{2\pi}{N} \mathbf{j} \cdot \mathbf{\xii}} \widehat{\uu}(\mathbf{\xii})
\end{equation}
Similarly, we define the discrete Fourier transform on the boundary grid as follows
\begin{equation}\label{fourier_transform_2}
    \widetilde{\F}(m_1,m_2) = \sum_{ (k_1,k_2) \in Z_{NP}^2} e^{-i \frac{2\pi}{NP} (k_1 m_1 + k_2 m_2 )} \F_{k_1,k_2}
\end{equation}
\begin{equation}\label{inv_Fourier_2}
    \F_{k_1,k_2}= \frac{1}{N^2 P^2} \sum_{(m_1,m_2) \in Z_{NP}^2} e^{i \frac{2\pi}{NP} (k_1 m_1 + k_2 m_2 )} \widetilde{\F}(m_1,m_2)
\end{equation}
The difference operators $\nabla_h$ and $\Delta_h$ become multiplication operators in Fourier space 
\begin{equation}
    (\widehat{\nabla}_h)_\alpha (\xii) = \frac{i}{h} \sin(\frac{2\pi}{N} \xi_\alpha), \quad \alpha = 1,2,3
\end{equation}
\begin{equation}
    (\widehat{\Delta}_h) (\xii) = -\frac{4}{h^2} \sum_{\alpha=1}^3 \sin^2(\frac{\pi}{N} \xi_\alpha)
\end{equation}
The following definition of the Fourier series representation of the
 function $\phi$ will turn out to be useful:
\begin{equation}\label{fourier_series}
    \phi(r) = \sum_{q = -\infty}^{\infty} \Phi(q) e^{i \frac{2\pi}{N} q r}
\end{equation}
\begin{equation}\label{fourier_coefficients}
    \Phi(q) = \frac{1}{N} \int_{-\frac{N}{2}}^{\frac{N}{2}} \phi(r) e^{-i \frac{2\pi}{N} q r} dr = \frac{1}{N} \int_{-2}^{2} \phi(r) e^{-i \frac{2\pi}{N} q r} dr
\end{equation}
Here, we reinterpret the IB 4-point delta function $\phi(r)$ as a Fourier series with period $N$, since $r$ is in units of meshwidth of the fluid grid, and in those units our periodic domain has length $N$ in each coordinate direction. Since the IB 4-point delta function has support $[-2,2]$, we should here assume that $N\geq 4$ to make the periodic extension possible.
By using equation \eqref{fourier_series} to rewrite equation \eqref{big_U} with $t = n\Delta t$
 in terms of $\Phi$,
 and by using the evenness of the function $\phi$,
 we get the following:
\begin{equation}\label{big_U_transform}
    \mathbf{U}^n_{k_1,k_2} = \sum_{q_1,q_2,q_3 = -\infty}^{\infty} \Phi(q_1) \Phi(q_2) \Phi(q_3) e^{i\frac{2\pi}{NP}(k_1 q_1 + k_2 q_2)} e^{i\frac{2\pi}{N}(\mathbf{q} \cdot \bm\epsilon)} \widehat{\mathbf{u}}^n (\mathbf{q})
\end{equation}
Now we multiply both sides of \eqref{big_U_transform} by $e^{-i\frac{2\pi}{NP}(k_1 m_1 + k_2 m_2)}$ and sum over $(k_1,k_2) \in Z_{NP}^2$. To evaluate this sum, we make use of 
\begin{equation}
    \sum_{k\in Z_{NP}} e^{i\frac{2\pi}{NP}k (q-m)} = \begin{dcases}
    NP, \hspace{5pt} \text{if} \hspace{5pt} q-m \hspace{5pt} \text{is an integer multiple of} \hspace{5pt} NP\\
    0, \hspace{5pt} \text{otherwise}\end{dcases}
\end{equation}
In this way, we get the result 
\begin{align*}
    \widetilde{ \mathbf{U} }^n (m_1,m_2) = & (NP)^2 e^{i\frac{2\pi}{N} (m_1 \epsilon_1 + m_2 \epsilon_2)} (\sum_{l_1 = -\infty}^{\infty} \Phi(m_1 + l_1 NP) e^{i 2\pi P l_1 \epsilon_1})\\ &  (\sum_{l_2 = -\infty}^{\infty} \Phi(m_2 + l_2 NP) e^{i 2\pi P l_2 \epsilon_2}) (\sum_{q_3 = -\infty}^{\infty} \Phi(q_3) e^{i\frac{2\pi}{N} q_3 \epsilon_3} \widehat{\mathbf{u}}^n (m_1,m_2,q_3)) \numberthis
\end{align*}
Similarly, we make use of \eqref{fourier_series} on the right-hand side of \eqref{leapfrog_2} , apply the definition \eqref{fourier_transform_2}, multiply both sides of the equation by $e^{-i\frac{2\pi}{N}(\xii \cdot \mathbf{j})}$, and sum over $\mathbf{j} \in Z_{N}^3$, we get the result
\begin{align*}\label{small_f_transform}
    & \widehat{\f}^{n+\frac{1}{2}} (\xii) = \frac{h_B^2 N^3}{h^3} e^{-i\frac{2\pi}{N}(\xii \cdot \bm \epsilon)}  \left(\sum_{l_3 = -\infty}^\infty \Phi(\xi_3 + l_3 N) e^{-i2\pi l_3 \epsilon_3}\right) \\ & \left(\sum_{l_1,l_2 = -\infty}^\infty \Phi(\xi_1 + l_1 N)\Phi(\xi_2 + l_2 N)e^{-i 2\pi (l_1 \epsilon_1 + l_2 \epsilon_2)} \widetilde{F}^{n+\frac{1}{2}}(\xi_1 + l_1 N,\xi_2 + l_2 N) \right) \numberthis
\end{align*}
Here we have used
\begin{equation}
    \sum_{k\in Z_{N}} e^{i\frac{2\pi}{N}k (q-\xi)} = \begin{dcases}
    N, \hspace{5pt} \text{if} \hspace{5pt} q-\xi \hspace{5pt} \text{is an integer multiple of} \hspace{5pt} N\\
    0, \hspace{5pt} \text{otherwise}\end{dcases}
\end{equation}
In equation \eqref{big_U_transform}, we can use the periodicity of $\widehat{\mathbf{u}}^n$ to rewrite the last factor as follows:
\begin{equation}
    \left( \sum_{m_3\in Z_N} (\sum_{l_3 = -\infty}^\infty \Phi(m_3+l_3 N) e^{i 2\pi l_3 \epsilon_3}) e^{i\frac{2\pi}{N}m_3 \epsilon_3} \widehat{\mathbf{u}}^n (\mathbf{m})        \right)
\end{equation}
where we let 
\begin{equation}
    q_3 = m_3 + l_3 N
\end{equation}
In a similar way, we can use the periodicity of $\widetilde{\F}$ to rewrite the sum over $(l_1,l_2)$ in equation \eqref{small_f_transform}. Note that $\widetilde{\F}$ is periodic in each of its arguments with period $NP$. Let 
\begin{equation}
    l_1 = p_1 + P l_1^{'}
\end{equation}
\begin{equation}
    l_2 = p_2 + P l_2^{'}
\end{equation}
Then, the sum over $(l_1,l_2)$ in equation \eqref{small_f_transform} can be rewritten as follows
\begin{align*}
    \sum_{(p_1,p_2)\in Z_P^2} & \left(\sum_{l_1 = -\infty}^\infty \Phi(\xi_1+NP_1 + NP l_1^{'}) e^{-i2\pi P l_1^{'} \epsilon_1}\right) \\ & \left(\sum_{l_2 = -\infty}^\infty \Phi(\xi_2+NP_2 + NP l_2^{'}) e^{-i 2\pi P l_2^{'} \epsilon_2}\right)\\ & e^{-2\pi (p_1 \epsilon_1 + p_2 \epsilon_2) \widetilde{\F}^{n+\frac{1}{2}} (\xi_1 + N p_1, \xi_2 + N p_2 )} \numberthis
\end{align*}
To simplify the notation, we let 
\begin{equation}\label{a}
    a(m,\epsilon) = \sum_{l = -\infty}^\infty \Phi(m+l NP) e^{i2\pi P l \epsilon}
\end{equation}
\begin{equation}\label{b}
    b(\xi,\epsilon) = \sum_{l = -\infty}^\infty \Phi(\xi+l N) e^{i2\pi l \epsilon}
\end{equation}
Since $\phi$ is real and even, $\Phi$ is real and even, we also have
\begin{equation}
    \overline{a(m,\epsilon)} = \sum_{l = -\infty}^\infty \Phi(m+l NP) e^{-i2\pi P l \epsilon}
\end{equation}
\begin{equation}
    \overline{b(\xi,\epsilon)} = \sum_{l = -\infty}^\infty \Phi(\xi+l N) e^{-i2\pi l \epsilon}
\end{equation}
Moreover, we let 
\begin{align*}
    \widehat{\F}^{n+\frac{1}{2}} (\xi_1,\xi_2,\epsilon_2,\epsilon_2) = \sum_{(p_1,p_2)\in Z_P^2} & \overline{a(\xi_1 + N p_1 ,\epsilon_1)} \hspace{3pt} \overline{a(\xi_2 + N p_2 ,\epsilon_2)} \\ & e^{-i2\pi (p_1 \epsilon_1 + p_2 \epsilon_2)} \widetilde{\F}^{n+\frac{1}{2}} (\xi_1 + N p_1 ,\xi_2 + N p_2) \numberthis
\end{align*}
Therefore, equation \eqref{big_U_transform} and equation \eqref{small_f_transform} become 
\begin{align*}
    \widetilde{\mathbf{U}}^n(\xi_1 + Np_1, \xi_2 + N p_2) = & (NP)^2 a(\xi_1 + N p_1,\epsilon_1) a(\xi_2 + N p_2,\epsilon_2) \\ & \left( \sum_{\xi_3 \in Z_N} b(\xi_3,\epsilon_3) e^{i\frac{2\pi}{N} (\bm \epsilon \cdot \xii)}  e^{i 2\pi (p_1 \epsilon_1 + p_2 \epsilon_2)} \widehat{\mathbf{u}}^n (\xii)\right) \numberthis
\end{align*}
\begin{equation}
    \widehat{\f}^{n+\frac{1}{2}} (\xii) = \frac{h_B^2 N^3}{h^3} e^{-i\frac{2\pi}{N} (\xii \cdot \bm \epsilon)}\widehat{\F}^{n+\frac{1}{2}} (\xi_1,\xi_2,\epsilon_2,\epsilon_2) \overline{b(\xi_3,\epsilon_3)}
\end{equation}
By equation \eqref{leapfrog_1}, we also have 
\begin{equation}
    \widetilde{\F}^{n+\frac{1}{2}}(\xi_1 + Np_1, \xi_2 + N p_2) = \widetilde{\F}^{n-\frac{1}{2}}(\xi_1 + Np_1, \xi_2 + N p_2) - \Delta t K \widetilde{\mathbf{U}}^n(\xi_1 + Np_1, \xi_2 + N p_2)
\end{equation}
After application of the discrete Fourier transform defined on the fluid grid, the system \eqref{leapfrog_3} becomes
\begin{equation}\label{sys_f}
\begin{dcases}
\rho \frac{\widehat{\uu}^{n+1} - \widehat{\uu}^n }{\dt} + \widehat{\nabla}_h \widehat{\mathbf{p}}^{n+\frac{1}{2}} = \mu \widehat{\Delta}_h \frac{\widehat{\uu}^{n+1} + \widehat{\uu}^n }{2} + \widehat{\mathbf{f}}^{n+\frac{1}{2}}\\
    \widehat\nabla_h \cdot \widehat{\mathbf{u}}^{n+1} = 0
\end{dcases}
\end{equation}
We can eliminate the pressure term $\widehat{\mathbf{p}}^{n+\frac{1}{2}}$ and solve for $\widehat{\uu}^{n+1}$ by introducing a $3\times3$ matrix 
\begin{equation}\label{before_proj}
    \widehat{P} (\mathbf{\xii}) = I - \frac{\left(\widehat\nabla_h(\mathbf{\xii})\right)  \left(\widehat\nabla_h(\mathbf{\xii})\right)^\ast}{\left(\widehat\nabla_h(\mathbf{\xii})\right)^\ast \left(\widehat\nabla_h(\mathbf{\xii})\right)}
\end{equation}
where $\ast$ denotes the Hermitian conjugate of a matrix and I is the $3\times 3$ identity matrix. The projection matrix $\widehat{P} (\mathbf{\xii})$ projects a vector field $\widehat{\uu}^{n+1}$ onto its divergence-free component by the Helmholtz decomposition. Therefore, $\widehat{P} (\mathbf{\xii}) \widehat{\nabla}_h \widehat{\mathbf{p}}^{n+\frac{1}{2}} = 0$ and $\widehat{P} (\mathbf{\xii}) \widehat{\uu}^{n+1} = \widehat{\uu}^{n+1}$. Then the system 
\eqref{before_proj} is reduced to 
\begin{equation}
    \rho \frac{\widehat{\uu}^{n+1} - \widehat{\uu}^n }{\dt}  = \mu \widehat{\Delta}_h \frac{\widehat{\uu}^{n+1} + \widehat{\uu}^n }{2} + \widehat{P} (\mathbf{\xii}) \widehat{\mathbf{f}}^{n+\frac{1}{2}}\\
\end{equation}
Solving for $\widehat{\uu}^{n+1}$ gives 
\begin{equation}
    \widehat{\uu}^{n+1} (\xii) = \frac{ \left(1 + \frac{\dt \mu}{2\rho} \widehat{\Delta}_h (\xii) \right) \widehat{\uu}^{n} (\xii) + \frac{\rho}{\dt} \widehat{P} (\mathbf{\xii}) \widehat{\mathbf{f}}^{n+\frac{1}{2}}}{1 - \frac{\dt \mu}{2\rho} \widehat{\Delta}_h (\xii)}
\end{equation}
To study the stability of the leapfrog scheme, we look for a solution in which all the variables are multiplied by a possibly complex number $z$ at each time step. The type of the solution that we seek is defined by 
\begin{equation}\label{solution_form_1}
    \widehat{\uu}^{n} (\xii) = z^n \widehat{\uu}^{0} (\xii)
\end{equation}
\begin{equation}\label{solution_form_2}
    \widetilde{\F}^{n+\frac{1}{2}} (m_1,m_2) = z^n \widetilde{\F}^{\frac{1}{2}} (m_1,m_2)
\end{equation}
\begin{equation}\label{solution_form_3}
    \widehat{\mathbf{f}}^{n+\frac{1}{2}} (\xii) = z^n \widehat{\mathbf{f}}^{\frac{1}{2}} (\xii) 
\end{equation}
Here the superscript on $z$ is actually a power, whereas all of other superscripts are merely labels indicating time steps. Then, after substituting (\ref{solution_form_1}-\ref{solution_form_3}) into the Fourier transform of the leapfrog scheme, we obtain the following system:
\begin{align*}\label{F_hat}
    \widehat{\F}^{\frac{1}{2}} (\xi_1,\xi_2,\epsilon_2,\epsilon_2) = \sum_{(p_1,p_2)\in Z_P^2} & \overline{a(\xi_1 + N p_1 ,\epsilon_1)} \hspace{3pt} \overline{a(\xi_2 + N p_2 ,\epsilon_2)} \\ & e^{-i2\pi (p_1 \epsilon_1 + p_2 \epsilon_2)} \widetilde{\F}^{\frac{1}{2}} (\xi_1 + N p_1 ,\xi_2 + N p_2) \numberthis
\end{align*}

\begin{equation}\label{target_force}
    \widetilde{\F}^{\frac{1}{2}} (\xi_1 + Np_1, \xi_2 + N p_2) =  \frac{-z}{z-1} \Delta t K \widetilde{\mathbf{U}}^0 (\xi_1 + Np_1, \xi_2 + N p_2)
\end{equation}

\begin{align*}\label{big_U_transform_final}
    \widetilde{\mathbf{U}}^0(\xi_1 + Np_1, \xi_2 + N p_2) = & (NP)^2 a(\xi_1 + N p_1,\epsilon_1) a(\xi_2 + N p_2,\epsilon_2) \\ & e^{i 2\pi (p_1 \epsilon_1 + p_2 \epsilon_2)} \left( \sum_{\xi_3 \in Z_N} b(\xi_3,\epsilon_3) e^{i\frac{2\pi}{N} (\bm \epsilon \cdot \xii)}   \widehat{\mathbf{u}}^0 (\xii)\right) \numberthis
\end{align*}

\begin{equation}\label{fluid_sys}
    \widehat{\uu}^{0} (\xii) = \frac{ \frac{\dt}{\rho} \widehat{P} (\mathbf{\xii}) \widehat{\mathbf{f}}^{\frac{1}{2}}(\xii)}{(z-1) - (z+1)\frac{\dt \mu}{2\rho} \widehat{\Delta}_h (\xii)}
\end{equation}

\begin{equation}\label{before_simp_force}
    \widehat{\f}^{\frac{1}{2}} (\xii) = \frac{h_B^2 N^3}{h^3} e^{-i\frac{2\pi}{N} (\xii \cdot \bm \epsilon)}\widehat{\F}^{\frac{1}{2}} (\xi_1,\xi_2,\epsilon_2,\epsilon_2) \overline{b(\xi_3,\epsilon_3)}
\end{equation}
Combining equation \eqref{F_hat} - \eqref{before_simp_force} gives
\begin{equation}\label{linear_sys}
    (\frac{(z-1)^2}{z} I + A) \widehat{F} ^{\frac{1}{2}} (\xi_1,\xi_2,\epsilon_2,\epsilon_2) = 0
\end{equation}
where 
\begin{align*}\label{matrix_A}
    A  =  \frac{N^5 \dt^2 K}{h \rho } & \left(\sum_{p_1\in Z_P}  \left|a(\xi_1 + N p_1 ,\epsilon_1)\right|^2\right) \left(\sum_{p_2\in Z_P} \left|a(\xi_2 + N p_2 ,\epsilon_2)\right|^2\right)\\&  \sum_{\xi_3 \in Z_N} \frac{ \widehat{P} (\mathbf{\xii})  \left|b(\xi_3,\epsilon_3)\right|^2}{1 - \frac{z+1}{z-1}\frac{\dt \mu}{2\rho} \widehat{\Delta}_h (\xii)}  \numberthis
\end{align*}
Thus, \eqref{linear_sys} has non-trivial solutions if and only if 
\begin{equation}\label{determinant}
    \det (\frac{(z-1)^2}{z} I + A) = 0
\end{equation}
In other words, $A$ must have an eigenvalue $\lambda = -\frac{(z-1)^2}{z}$. We expect that all solutions $z$ lie strictly inside the unit circle for $\dt$ positive and sufficiently small. The mechanism of instability then has to be some solution $z$ crossing the unit circle for some $(\xi_1,\xi_2)$. We claim that this can only happen when $z = -1$. To see this, note for $z = e^{i\theta}$ on the unit circle, we have 
\begin{equation}
    \frac{(z-1)^2}{z} = 2(\cos\theta - 1)
\end{equation}
which is real, but $\frac{z+1}{z-1} = \frac{\cos(\frac{\theta}{2})}{i\sin(\frac{\theta}{2})}$
is a nonzero imaginary number unless $\theta = \pi$, which is the same as $z = -1$. From the form of (\ref{matrix_A}), it therefore seems clear that $A$ cannot have real eigenvalues when $z$ is on the unit circle unless $z = -1$ and therefore that the only nontrivial way to satisfy \eqref{linear_sys} for z on the unit circle is if $z = -1$. When $z = -1$, the matrix $A$ becomes
\begin{align*}
    A_0 = &  \frac{N^5 \dt^2 K}{h \rho }  \left(\sum_{p_1\in Z_P}  \left|a(\xi_1 + N p_1 ,\epsilon_1)\right|^2\right) \left(\sum_{p_2\in Z_P} \left|a(\xi_2 + N p_2 ,\epsilon_2)\right|^2\right) \\ & \sum_{\xi_3\in Z_N} \left|b(\xi_3,\epsilon_3)\right|^2 \widehat{P}(\xii) \numberthis
\end{align*}
which is a linear combination of projection matrices, and moreover, equation \eqref{determinant} becomes simply the statement that $A_0$ has an eigenvalue equal to $4$. Note that when $z = -1$, the viscosity term drops out because the fluid solver treats the diffusion term implicitly through the trapezoidal rule. The $3\times 3$ matrix $A_0$ is real, symmetric, and nonnegative, and its eigenvalues are bounded by\begin{equation}\label{lambda_max}
    \lambda_{max} (A_0) \leq \frac{K(\Delta t)^2}{\rho h} C(\xi_1,\xi_2)
\end{equation}
where 
\begin{align*}
    C(\xi_1,\xi_2) = & N^5 \left(\sum_{p_1\in Z_P}  \left|a(\xi_1 + N p_1 ,\epsilon_1)\right|^2\right) \left(\sum_{p_2\in Z_P} \left|a(\xi_2 + N p_2 ,\epsilon_2)\right|^2\right)\\ & \sum_{\xi_3\in Z_N} \left|b(\xi_3,\epsilon_3)\right|^2
    \numberthis
\end{align*}
Now let $(\Delta_t)_c >0$ be defined by
\begin{equation}\label{max_C}
    \frac{K(\Delta_t)_c^2}{\rho h} \max_{(\xi_1,\xi_2) \in Z_N^2} C(\xi_1,\xi_2) = 4
\end{equation}
Then, $\Delta t \in \left(0, (\Delta t)_c \right)$ is a sufficient condition for stability, since all solutions $z$ of \eqref{determinant} lie inside the unit circle for $\Delta t$ positive and sufficiently small, and since they cannot escape from within the unit circle unless $A_0$ has an eigenvalue equal to 4, and since this is impossible for $\Delta t < (\Delta t)_c$ because of \eqref{lambda_max} - \eqref{max_C}.

\begin{figure}[h]
\centering
\includegraphics[width=0.9\textwidth]{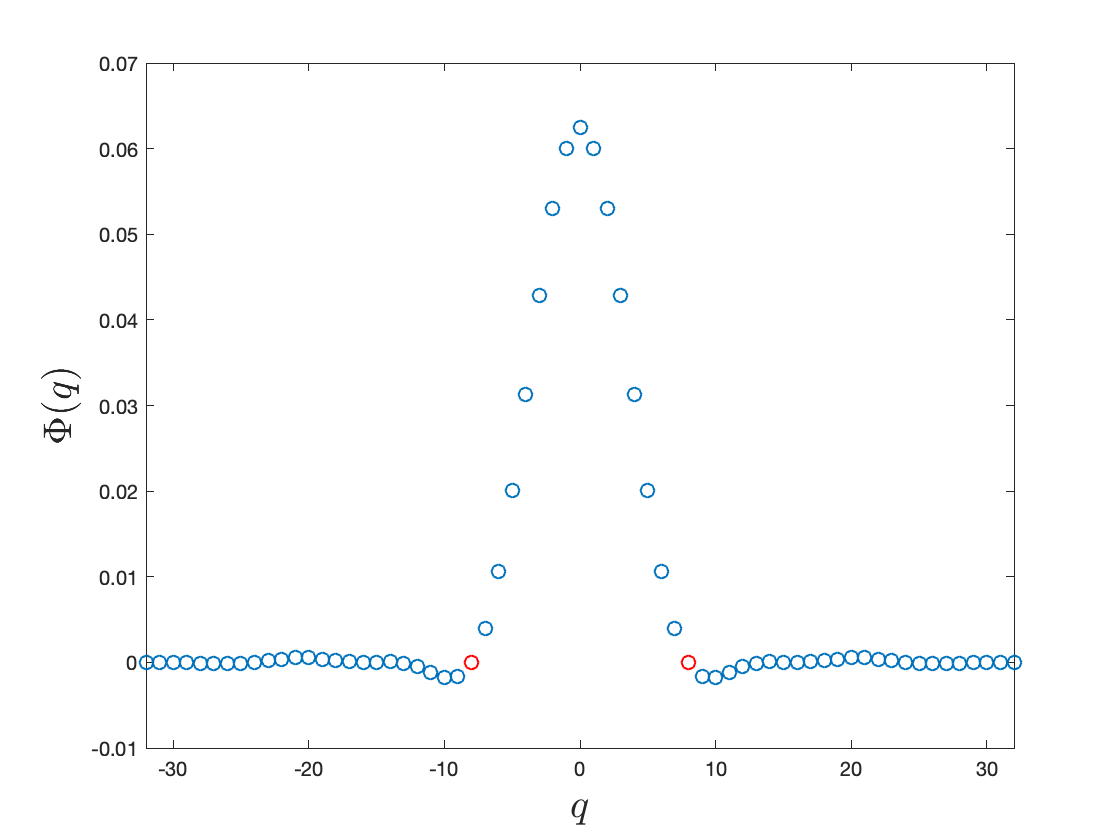}
\caption{A plot of $\Phi(q)$ for N = 16. The red circles correspond to $\Phi(\frac{N}{2})$ and $\Phi(-\frac{N}{2})$. From this plot, it is clear that $\Phi(q)$ bell-shaped and most of its mass lies within $(-\frac{N}{2},\frac{N}{2})$. Note that
$\Phi(q)$ is defined for all integer values of $q$; in the figure only
the values of $q \in [-2N,2N]$ are shown.}
\label{fig:1}
\end{figure}

To evaluate $\max_{(\xi_1,\xi_2) \in Z_N^2} C(\xi_1,\xi_2)$, we make the \textit{band-limited approximation} here, which states
\begin{equation}
    \Phi(q) = 0 \hspace{5pt} \text{for} \hspace{5pt} |q| >\frac{N}{2}
\end{equation}
The motivation of making this band-limited approximation is that $\Phi(q)$ is a bell-shaped function and is approximately band-limited. A plot of $\Phi(q)$ is given in Figure \ref{fig:1}.

\begin{figure}[h]
\centering
\includegraphics[width=0.9\textwidth]{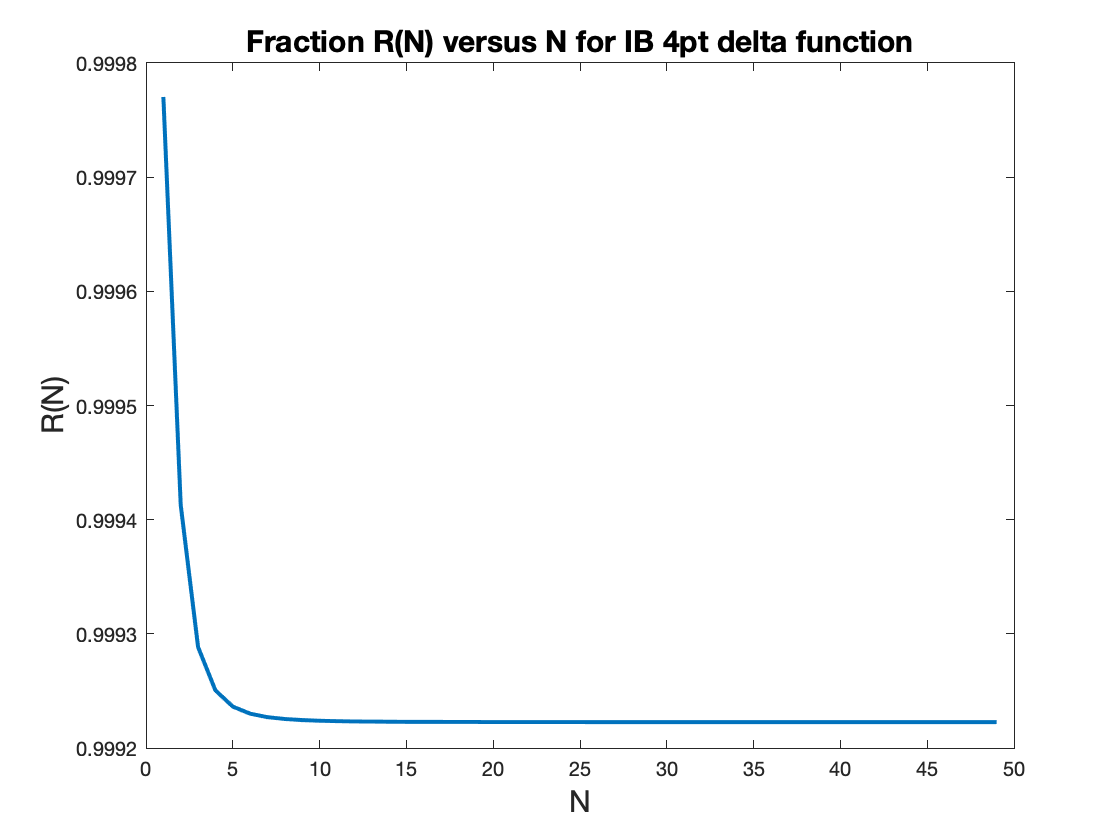}
\caption{A plot of $R(N)$ for $N \geq 4$. Note the scale of the vertical axis, which suggests that $R(N)$ is very close to 1 for all $N$'s.}
\label{fig:2}
\end{figure}

We also investigate how accurate the band-limited approximation is in a quantitative way by evaluating the following ratio:
\begin{equation}
    R(N) = \frac{\sum_{p=-\frac{N}{2}}^{\frac{N}{2}} \Phi^2(p) }{\sum_{q=-\infty}^{\infty} \Phi^2(q)}
\end{equation}
as a function of $N$. Note that as $N\rightarrow \infty$, the ratio $R(N)$ converges to a limit and $\Phi(-\frac{N}{2}) = \Phi(\frac{N}{2}) = 0$. Whether including the boundary terms in the numerator of $R(N)$ or not does not change the value of $R(N)$. Let $\widehat \phi$ denote the Fourier transform of the IB 4-point delta function $\phi$ and it is defined by
\begin{equation}
    \widehat{\phi}(s) = \int_{-\infty}^\infty \phi(r) e^{-isr} dr = \int_{-2}^2 \phi(r) e^{-isr} dr
\end{equation}
For any positive integer $N$, we have the following correspondence between the Fourier transform and the Fourier coefficients of $\phi(r)$:
\begin{equation}
    \Phi(q) = \frac{1}{N} \widehat\phi(\frac{2\pi}{N}q)
\end{equation}
Therefore, as $N \rightarrow \infty$, we have 
\begin{equation}
    R(N) \rightarrow \frac{\int_{-\pi}^\pi \widehat{\phi}^2(s) ds}{\int_{-\infty}^\infty \widehat{\phi}^2(s) ds}
\end{equation}
In practice, we do not need to compute the integral in the denominator of $R(N)$ because of the Parseval's theorem and the sum of squares property of the IB 4-point delta function\cite{peskin_2002}, which gives
\begin{equation}
    \sum_{q=-\infty}^{\infty} \Phi^2(q) = \frac{3}{8N}
\end{equation}
\begin{equation}
    \int_{-\infty}^\infty \widehat{\phi}^2(s) ds = \frac{3}{8}
\end{equation}
We compute $R(N)$ for different values of $N$, and a plot of it as a function of $N$ is given in Figure \ref{fig:2}, which numerically justifies the band-limited approximation of $\phi(r)$.
Let $A(m) = \sum_{l = -\infty}^\infty |\Phi(m+lNP)|$ and by the triangular inequality, we obtain
\begin{equation}
    |a(m,\epsilon)| \leq \sum_{l=-\infty}^\infty |\Phi(m+lNP)| = A(m)
\end{equation}
Since $A(m)$ is periodic with period $NP$, we can restrict consideration to $-\frac{NP}{2} \leq m \leq \frac{NP}{2}$ when evaluating $A(m)$. Then, for $l\neq 0$, we have $|m+ lNP| > \frac{NP}{2} \geq \frac{N}{2}$. So, by the band-limited approximation, $l\neq 0$ implies $\Phi(m+lNP) = 0$. It follows that
\begin{equation}
    A(m) = |\Phi(m)|, \hspace{5pt} -\frac{NP}{2} \leq m \leq \frac{NP}{2}
\end{equation}
Again, by the band-limited approximation, we get 
\begin{equation}
    \sum_{p\in Z_P}  \left|a(\xi + N p ,\epsilon)\right|^2 \leq \sum_{p\in Z_P}  A^2(\xi + N p) = A^2(\xi) = |\Phi(\xi)|^2, \hspace{5pt} -\frac{N}{2} \leq \xi \leq \frac{N}{2}
\end{equation}
Since $\phi(r)$ is real and positive, we have 
\begin{equation}
    |\Phi(\xi)| = \frac{1}{N} |\int_{-2}^{2} \phi(r) e^{-i \frac{2\pi}{N} \xi r} dr| \leq  \frac{1}{N} \int_{-2}^{2} |\phi(r)| dr= \frac{1}{N} = \Phi(0)
\end{equation}
Therefore, 
\begin{equation} \label{bound_for_a}
     \sum_{p\in Z_P}  \left|a(\xi + N p ,\epsilon)\right|^2 \leq \frac{1}{N^2}
\end{equation}
Next, by the definition of $b(\xi,\epsilon)$, we get
\begin{align*}\label{simp_b}
    b(\xi,\epsilon) & = \sum_{l=-\infty}^\infty \Phi(\xi + l N) e^{i2\pi l\epsilon} \\ & = \frac{1}{N} \sum_{l=-\infty}^\infty \left(\int_{-2}^{2} \phi(r) e^{-i\frac{2\pi}{N} (\xi + lN) r} dr\right) e^{i2\pi l\epsilon} \\ & = \frac{1}{N}  \int_{-2}^{2} \phi(r) e^{-i\frac{2\pi}{N} \xi r} \sum_{l=-\infty}^\infty e^{i2\pi l(\epsilon-r)} dr \\ & = \frac{1}{N}  \int_{-2}^{2} \phi(r) e^{-i\frac{2\pi}{N} \xi r} \sum_{l=-\infty}^\infty \delta(j+\epsilon-r) dr  \numberthis
\end{align*}
Here $\delta$ represents the Dirac delta function and the following identity was used in \eqref{simp_b}
\begin{equation}\label{dirac_delta_fourier}
    \sum_{j=-\infty}^\infty \delta(j+x) = \sum_{l=-\infty}^\infty e^{i2\pi lx}
\end{equation}
The left-hand side of \eqref{dirac_delta_fourier} is a periodic series with period $1$, and the right-hand side of \eqref{dirac_delta_fourier} is its Fourier series. Therefore, it follows that
\begin{equation}
    b(\xi,\epsilon) = \frac{1}{N} \sum_{j:j+\epsilon\in(-2,2)} \phi(j+\epsilon) e^{-i\frac{2\pi}{N}\xi(j+\epsilon)}
\end{equation}
Note that we can write the interval as $(-2,2)$ because $\phi(2) = \phi(-2) = 0$.
It follows that
\begin{equation}\label{sum_for_b}
    |b(\xi,\epsilon)|^2 = \frac{1}{N^2} \sum_{j:j+\epsilon\in(-2,2)} \sum_{k:k+\epsilon\in(-2,2)} \phi(j+\epsilon) \phi(k+\epsilon) e^{-i\frac{2\pi}{N} \xi (j-k)}
\end{equation}
Note that $\epsilon$ only appears in the argument of $\phi$; it has canceled out of the exponent. Moreover, \eqref{sum_for_b} is exact and we do not need to apply the band-limited approximation here. Now sum over $\xi \in Z_N$ and note that 
\begin{equation}
    \sum_{\xi\in Z_N} e^{-i\frac{2\pi}{N}\xi(j-k)} = \begin{dcases}
    N, \hspace{5pt} \text{if} \hspace{5pt} j-k \hspace{5pt} \text{is an integer multiple of} \hspace{5pt} N\\
    0, \hspace{5pt} \text{otherwise}\end{dcases}
\end{equation}
Given the restriction that $j+\epsilon \in (-2,2)$ and $k+\epsilon \in (-2,2)$, we have $|j-k| < 4$. So, if $N \geq 4$, the only way that $j-k$ can be a multiple of $N$ is if $j-k=0$. By the sum of squares property of the IB 4-point delta function\cite{peskin_2002}, it follows that
\begin{equation}\label{bound_for_b}
    \sum_{\xi\in Z_N} |b(\xi,\epsilon)|^2 = \frac{1}{N} \sum_{j:j+\epsilon \in (-2,2)} \phi^2(j+\epsilon) = \frac{1}{N}(\frac{3}{8}) = \frac{3}{8N}
\end{equation}
Thus, by the two bounds we have obtained and stated in \eqref{bound_for_a} and \eqref{bound_for_b}, 
\begin{equation}
    \max_{(\xi_1,\xi_2) \in Z_N^2} C(\xi_1,\xi_2) = C(0,0) \leq \frac{3}{8}
\end{equation}
The stability criterion \eqref{max_C} can be rewritten as an \textit{extremely} simple expression
\begin{equation}\label{stability_criteria_1}
    \frac{ \dt^2 K}{ \rho h}  \leq \frac{32}{3}
\end{equation}
We will later show in Section 3 that the stability criterion does not only works for the time-dependent Stokes equations we are considering but also for the Navier-Stokes equations, which are used in most of the simulations. 

Note that a direct consequence of \eqref{stability_criteria_1} is that we can achieve the continuum limit ($h \rightarrow 0$, $\Delta t\rightarrow 0$) and the no-slip limit ($K\rightarrow \infty$) simultaneously by letting $\Delta t \rightarrow 0, K\propto \frac{1}{\Delta t}$, and $h\propto \Delta t$. The stability criterion \eqref{stability_criteria_1} ensures numerical stability of doing so. 

Also note that the band-limited approximation is only applied when getting an estimate of $|a(m,\epsilon)|$. For the speical case $P = 1$, $a(m,\epsilon) = b(m,\epsilon)$ for any $m$ and $\epsilon$. So, the band-limited approximation is no longer needed in the case $P = 1$ and for that case, the stability criterion \eqref{stability_criteria_1} is exact. 

\subsection{Generalization to a planar elastic membrane}
Now, with the same configuration and the same numerical scheme, we present a generalization of our analysis framework to the case of a planar elastic membrane. The equations of the motion are the same as the ones in Section 2.1 except for the equation of the force on the immersed boundary, which is stated as below: 
\begin{equation}\label{force_eq_membrane}
    \mathbf{F} (x_1,x_2,t) = K \Delta X (s_1,s_2,t)
\end{equation}
where $\Delta$ denotes the Laplacian operator. The spatial discretization of equation \eqref{force_eq_membrane} is defined as follows
\begin{align*}\label{discretized_force_membrane}
    \mathbf{F}_{k_1,k_2} (t) & = K \left(\frac{\X_{k_1+1,k_2}(t) + \X_{k_1-1,k_2}(t) - 2\X_{k_1,k_2}(t)}{h_B^2}   \right) \\ & + K \left(\frac{\X_{k_1,k_2+1}(t) + \X_{k_1,k_2-1}(t) - 2\X_{k_1,k_2}(t)}{h_B^2}   \right) \numberthis
\end{align*}
Note that we can no longer evaluate the IB 4-point delta function $\delta_h$ at fixed positions $\mathbf{X}^0 (k_1,k_2)$ in practice because the boundary is elastic and motion of the boundary is allowed. But still, we evaluate the IB 4-point delta function $\delta_h$ at fixed positions $\mathbf{X}^0 (k_1,k_2)$ here for the stability analysis. This means we are only considering small-amplitude vibrations of an elastic membrane. In Section 3.2, we will show that in this way we get a good approximation to the actual stability behavior observed numerically, even for vibrations of larger amplitudes. Since we are still evaluating the IB 4-point delta function at fixed positions and we have here dropped the convective term of the Naiver-Stokes equations, the Runge-Kutta scheme reduces to a simple leapfrog scheme, which is the same as what is described in Section 2.2 except for the equation of the boundary force. 

Using the same Fourier techniques as in Section 2.3, we can get the following result:
\begin{equation}\label{linear_sys_membrane}
    (\frac{(z-1)^2}{z} I + A) \widehat{F} ^{\frac{1}{2}} (\xi_1,\xi_2,\epsilon_1,\epsilon_2) = 0
\end{equation}
for some matrix $A$. We claim that $A$ cannot have real eigenvalues when $z$ is on the unit circle unless $z = -1$. For $z = -1$, we get an expression for $A$ as follows
\begin{align*}\label{matrix_A_membrane}
    A = &  \frac{4 N^5 \dt^2 K P^2}{ \rho h^3}  \sum_{p_1\in Z_P}  \left|a(\xi_1 + N p_1 ,\epsilon_1)\right|^2 \sum_{p_2\in Z_P} \left|a(\xi_2 + N p_2 ,\epsilon_2)\right|^2 \\ & \left(\sin^2(\frac{\pi}{NP} (\xi_1 + N p_1)) + \sin^2(\frac{\pi}{NP} (\xi_2 + N p_2) )\right ) \left(\sum_{\xi_3\in Z_N} \left|b(\xi_3,\epsilon_3)\right|^2 \widehat{P}(\xii)\right)   \numberthis
\end{align*}
Now, we let 
\begin{align*}
    C(\xi_1,\xi_2) = &  N^5 \sum_{p_1\in Z_P}  \left|a(\xi_1 + N p_1 ,\epsilon_1)\right|^2 \sum_{p_2\in Z_P} \left|a(\xi_2 + N p_2 ,\epsilon_2)\right|^2 \\ & \left(\sin^2(\frac{\pi}{NP} (\xi_1 + N p_1)) + \sin^2(\frac{\pi}{NP} (\xi_2 + N p_2) )\right ) \left(\sum_{\xi_3\in Z_N} \left|b(\xi_3,\epsilon_3)\right|^2 \right)   \numberthis
\end{align*}
Let the critical time step $(\Delta t)_c$ be defined by
\begin{equation}\label{stability_criteria_2}
    \frac{4 P^2 }{h^2}\frac{K (\Delta t)_c^2}{\rho h} \max_{(\xi_1,\xi_2) \in Z_N^2} C(\xi_1,\xi_2) = 4
\end{equation}
Then, $\Delta t \in (0, (\Delta t)_c)$ becomes a sufficient condition for stability. However, the maximizers of $C(\xi_1,\xi_2)$ will not be the same as the ones we found in Section 2.3 (i.e. $\xi_1 = \xi_2 = 0$) because of the sine terms. Note that if we absorb the term $\frac{4 P^2}{ h^2}$ into $C(\xi_1,\xi_2)$ and compare \eqref{stability_criteria_2} to \eqref{stability_criteria_1}, we can observe that the only difference between these equations is the term 
\begin{equation}
    \frac{4 P^2 }{h^2} \left(\sin^2(\frac{\pi}{NP} (\xi_1 + N p_1)) + \sin^2(\frac{\pi}{NP} (\xi_2 + N p_2) )\right )
\end{equation}
which is from the Fourier transform of equation \eqref{discretized_force_membrane}. 

Maximizing $C(\xi_1,\xi_2)$ over $(\xi_1,\xi_2) \in Z_N^2$ is hard because there is a trade-off between maximizing the sine terms and the sums that involve $a$ and $b$. Optimizing $C(\xi_1,\xi_2)$ numerically turns out to be a much simpler task because $Z_N^2$ is a finite set of numbers and we can compute $C(\xi_1,\xi_2)$ for each $(\xi_1,\xi_2) \in Z_N^2$. However, recalling the definitions of $a(m,\epsilon)$ given in equation \eqref{a} and \eqref{b}, we find that the numerical evaluation of $C(\xi_1,\xi_2)$ includes summing over all integers. To make the computation more efficient, we apply the band-limited approximation introduced in Section 2.3, and then $C(\xi_1,\xi_2)$ can be simplified as 
\begin{equation}
    C(\xi_1,\xi_2) = \frac{3 N^4 }{8}  |\Phi(\xi_1)|^2 |\Phi(\xi_2)|^2 \left(\sin^2(\frac{\pi}{NP} \xi_1) + \sin^2(\frac{\pi}{NP} \xi_2 )\right ) 
\end{equation}
where the bound $\left(\sum_{\xi_3\in Z_N} \left|b(\xi_3,\epsilon_3)\right|^2 \right) \leq \frac{3}{8N}$ is given by \eqref{bound_for_b}.  
Table \ref{table:wavenumber} gives the maximum values of $C(\xi_1,\xi_2)$ for different values of $N$ and $P$ and their maximizers (i.e. the most unstable wavenumbers). It is shown in Table \ref{table:wavenumber} that the maximum values of $C$ are almost invariant to $N$ but they depend on the value of $P$.

In this case, we can get to the continuum limit by letting $\Delta t\rightarrow 0$ and $h\propto \Delta t^{\frac{2}{3}}$ while still remaining stable. Note that the difference in this case of an elastic membrane and the case of no-slip boundary we analyzed in Section 2.3 is that the constant $K$ has a different unit and a different meaning. Here, $K$ is a physical parameter and we should not push it to infinity. Although the fact that $h\propto \Delta t^{\frac{2}{3}}$ makes the numerical scheme more unstable than the one analyzed in Section 2.3, it is easier to achieve the desired continuum limit in this case because $K$ does not need to be changed.
\begin{table*}[t!]
\centering
\begin{tabular}{| c| c| c|c|}
\hline
$N$ & $P$ & $\max C(\xi_1,\xi_2)$ & $\text{the most unstable wavenumbers}$\\ \hline
$16$ & $1$ & $5.786 \times 10^{-2}$ & $\xi_1 = 2, \xi_2 = 3$\\ \hline
$32$ & $1$ & $5.939 \times 10^{-2}$& $\xi_1 = 5, \xi_2 = 5$\\ \hline
$64$ & $1$ & $5.969 \times 10^{-2}$& $\xi_1 = 9, \xi_2 = 10$\\ \hline
$128$ & $1$ & $1.555 \times 10^{-2}$& $\xi_1 = 19, \xi_2 = 20$ \\ \hline
$16$ & $2$ & $1.555 \times 10^{-2}$& $\xi_1 = 2, \xi_2 = 3$\\ \hline
$32$ & $2$ & $1.574 \times 10^{-2}$& $\xi_1 = 4, \xi_2 = 5$\\ \hline
$64$ & $2$ & $1.578 \times 10^{-2}$& $\xi_1 = 10, \xi_2 = 11$\\ \hline
$128$ & $2$ & $1.578 \times 10^{-2}$& $\xi_1 = 19, \xi_2 = 20$\\ \hline
$16$ & $3$ & $7.005 \times 10^{-3}$& $\xi_1 = 2, \xi_2 = 3$\\ \hline
$32$ & $3$ & $7.074 \times 10^{-3}$& $\xi_1 = 4, \xi_2 = 5$\\ \hline
$64$ & $3$ & $7.084 \times 10^{-3}$& $\xi_1 = 10, \xi_2 = 11$\\ \hline
$128$ & $3$ & $7.089 \times 10^{-3}$& $\xi_1 = 19, \xi_2 = 20$\\ \hline
\end{tabular}
\caption{$\max C(\xi_1,\xi_2)$ for different values of $N$ and $P$ }\label{table:wavenumber}
\end{table*}
\section{Numerical results}

\subsection{Target point force}

We test the numerical stability analysis in Section 2.3 with fixed viscosity $\mu$ and density $\rho$. We start from a random velocity field drawn i.i.d. from a standard Gaussian distribution and average the results over 10 tests.
\begin{figure}[H]
\centering
\includegraphics[width=0.65\textwidth]{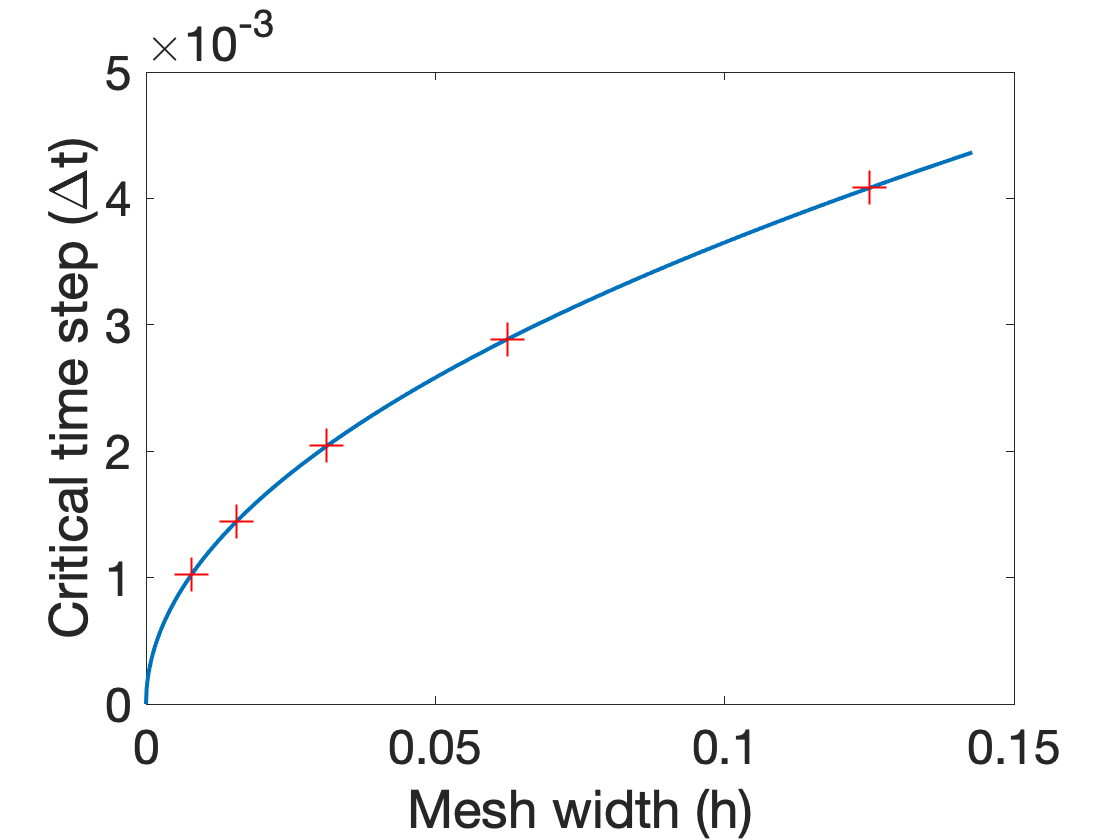}
\caption{Critical time steps for different mesh widths with $K = 8\times 10^4$ and $\rho$ = 1. The blue curve is what we would expect based on \eqref{stability_criteria_1} and red markers are data points we get from numerical tests.}
\label{fig:4}
\end{figure}
Figure \ref{fig:4} shows that the critical time step we observe numerically perfectly agrees with what the theory implies if we fix all the other parameters and vary the mesh width $h$. Moreover, for the parameter regime plotted on Figure \ref{fig:4}, the absolute difference of the numerically observed critical time step and the theoretical prediction is below 6 digits. We test different time steps and search for the critical time step by applying the bisection method. To determine whether the scheme is stable for a certain time step, we plot the ratio of the $L^2$ norm of the velocity at all time steps and the initial $L^2$ norm of the velocity for a sufficiently long period of time. Since there is no driving force, the energy (i.e. $L^2$ norm) of the flow should monotonically decay over time if the scheme is stable. If we enter the instability region, the energy (i.e. $L^2$ norm) of the flow should blow up. However, if we are in the instability region but are very close to the stability boundary, we may observe that the energy may decay transiently and then blow up. Here, we give a plot of the $L^2$ norm of the flow versus time for each pattern of growth/decay. 
\begin{figure}[H]
\centering
\includegraphics[width=0.65\textwidth]{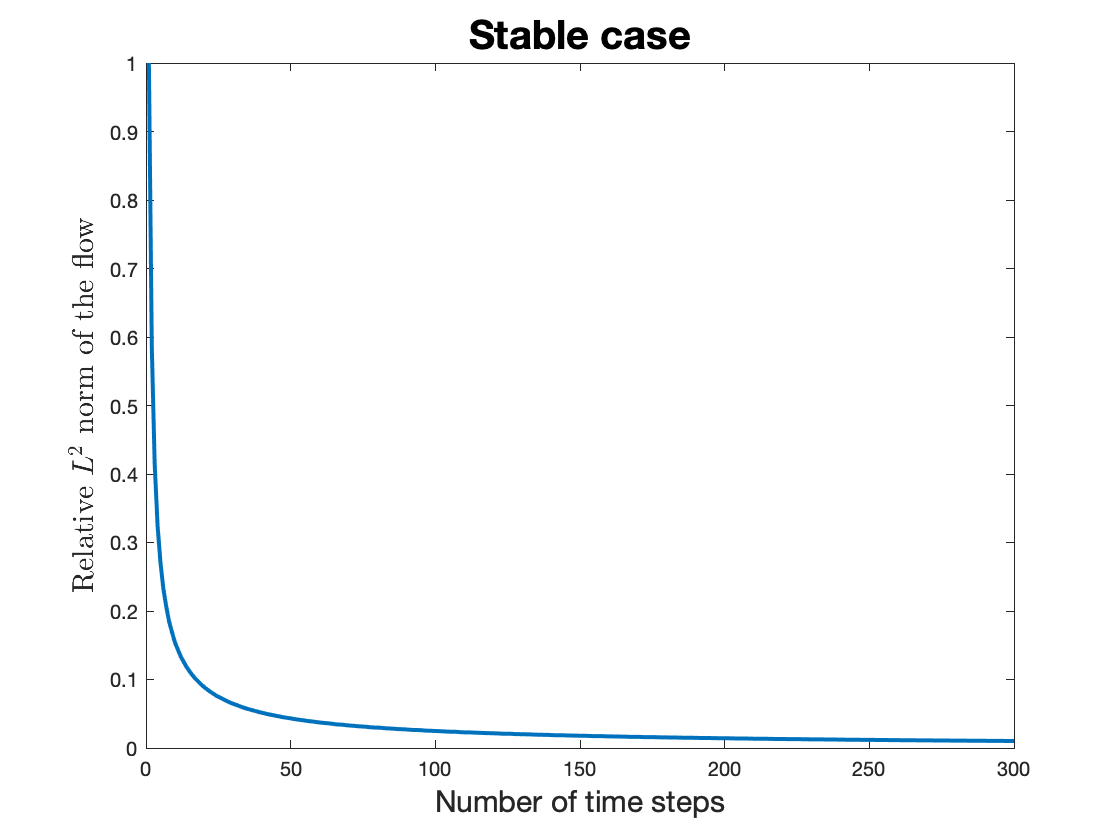}
\caption{Evolution of the relative $L^2$ norm ($L^2$ norm of the flow divided by the $L^2$ norm of the initial velocity) when the time step is small enough such that the numerical scheme is stable. Parameter values: $K = 8 \times 10^{4}$, $N = 32$, $P = 2$, $\Delta t = 2.0410 \times 10^{-3}$. Predicted critical time step: $\Delta t_c = 2.0412 \times 10^{-3}$. }
\label{fig:9}
\end{figure}

\begin{figure}[h]
\centering
\includegraphics[width=0.65\textwidth]{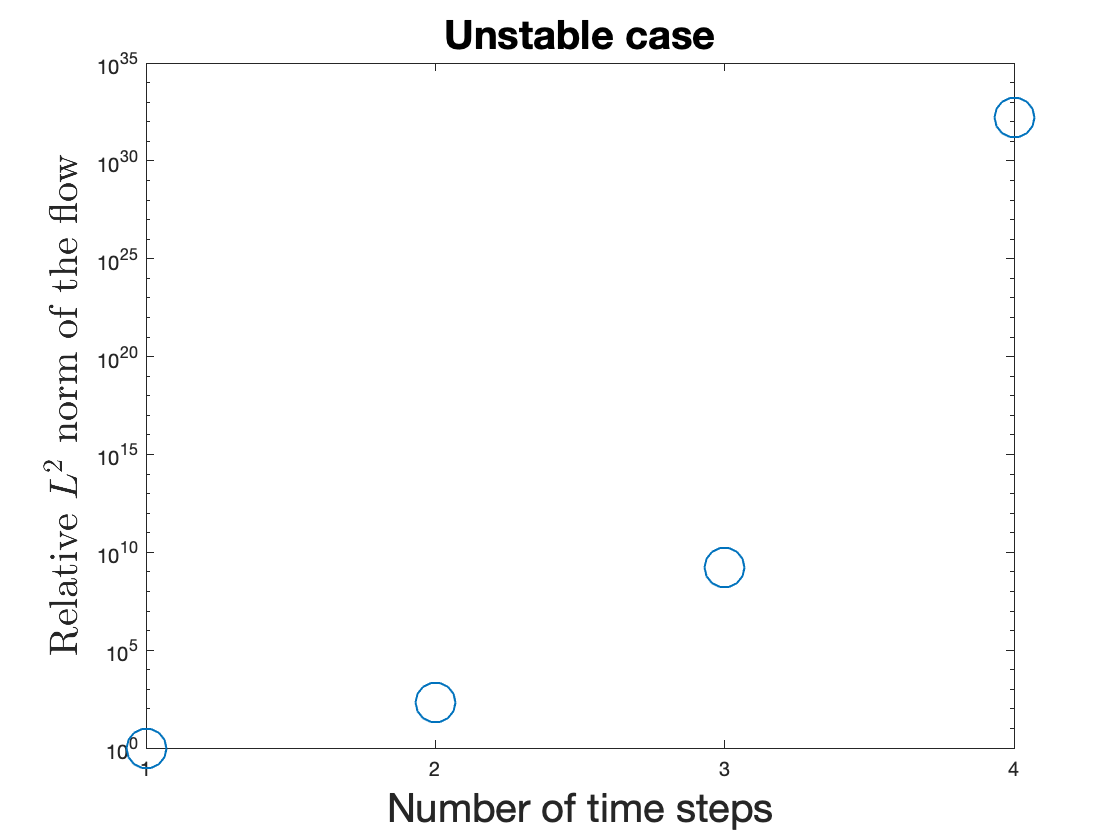}
\caption{Evolution of the relative $L^2$ norm ($L^2$ norm of the flow divided by the $L^2$ norm of the initial velocity) when the time step is large enough that the numerical scheme is unstable. The relative $L^2$ norm blows up in four time steps. Note that here we use a log plot on the y-axis. Parameter values: $K = 8 \times 10^{4}$, $N = 32$, $P = 2$, $\Delta t = 2.5410 \times 10^{-3}$. Predicted critical time step: $\Delta t_c = 2.0412 \times 10^{-3}$.  }
\label{fig:10}
\end{figure}

\begin{figure}[H]
\centering
\includegraphics[width=0.65\textwidth]{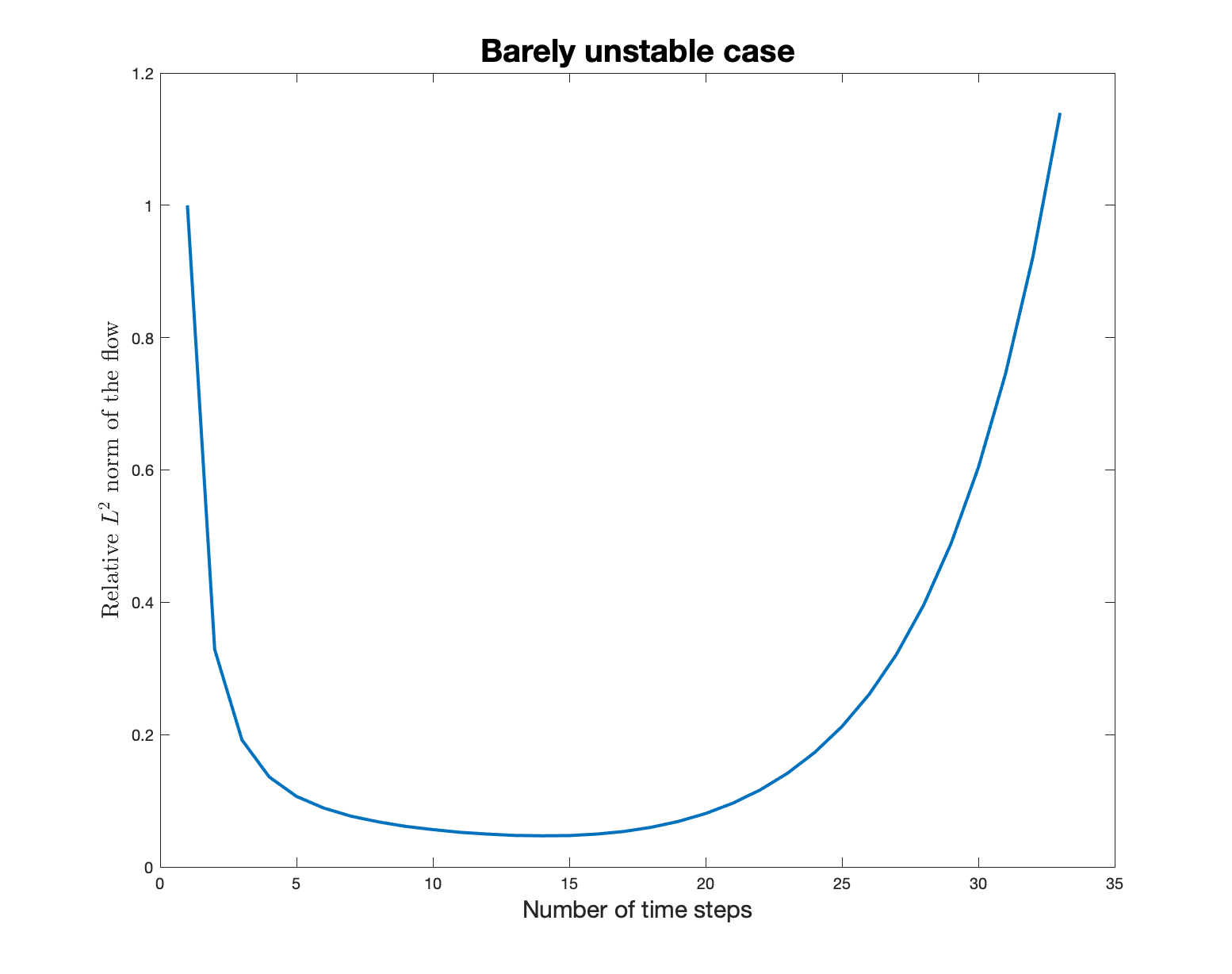}
\caption{Evolution of the relative $L^2$ norm ($L^2$ norm of the flow divided by the $L^2$ norm of the initial velocity) when the time step is close to the stability boundary but not enough small such that the scheme is stable. Parameter values: $K = 8 \times 10^{4}$, $N = 32$, $P = 2$, $\Delta t = 2.1430 \times 10^{-3}$. Predicted critical time step: $\Delta t_c = 2.0412 \times 10^{-3}$.  }
\label{fig:11}
\end{figure}

Moreover, for this parameter regime, if we replace the Stokes equations by the Navier-Stokes equations in numerical tests, the results do not change much. The numerically observed critical time step and the theoretical prediction still agree up to 5 digits. 

Note that \eqref{stability_criteria_1} implies a very important fact that the stability boundary does not depend on the density of the target points $P$. In practice, to get better volume conservation and avoid leakage of fluid, we need to place a large number of the target points on the no-slip boundary. Therefore, if the numerical result also suggests that the stability boundary is invariant to $P$, then we do not need to worry about numerical instability when we choose $P$. The only obstacle of increasing $P$ would be the computational cost for spreading the boundary force on each target point. In numerical tests, we take $P = 1,2,3,\cdots,10$ and the numerical stability boundary does not change, which verifies the amazing fact that the stability boundary does not depend on P. 

As the results about the stability boundary also apply to the Naiver-Stokes equations, we conclude that analyzing the time-dependent Stokes equations can give an extremely well approximation to the stability boundary of the Navier-Stokes equations.

\subsection{Achieving the Continuum Limit and the No-Slip Limit}

In this section, we provide a numerical example of a 3D Poiseuille flow in a periodic cube to demonstrate the fact that we can achieve the continuum limit ($h \rightarrow 0$, $\Delta t\rightarrow 0$) and the no-slip limit ($K\rightarrow \infty$) simultaneously by letting $\Delta t \rightarrow 0, K\propto \frac{1}{\Delta t}$, and $h\propto \Delta t$. 

In the numerical example, a 3D Poiseuille flow is simulated in a cubic box with periodic boundary conditions in all three directions and with the no-slip boundary condition implemented on the top and the bottom of the domain (i.e. $z = 0$ and $z = L$). A driving force is applied in the x-direction and it is uniform in both space and time. Therefore, a Poiseuille flow $\mathbf u(x,y,z) = \mathbf u(x,z)$ can be obtained and we know the analytical solution of it. 

We run the simulation on different grids. To keep staying in the stability region given by the analysis in Section 2.3, whenever we refine the grid by a factor of $2$, we increase $K$ by a factor of $2$ and decrease $\dt$ by a factor of $2$. The domain we use for all grids is $[0,1] \times [0,1] \times[0,1]$ and the coarsest grid size is $16\times 16 \times 16$ with $h_0 = 1/16$, $\Delta t_0 = 1/400$, and $K_0 = 8 \times 10^{4}$. We also keep $P = 2$ fixed. To measure the error on the boundary condition, we define the following metric:
\begin{equation}
    \mathbf{d} (\mathbf{X},\mathbf{Z}) = \mathbf{X}-\mathbf{Z}
\end{equation}
which is the displacement of the immersed boundary points $\X$ at the final time $t = 40$ from their initial and target position $\Z$. For simulations on all grids, we set the simulation time be $40$ and we start with steady-state solutions. The simulation time $t=40$ is enough for viscous effects to diffuse across the domain (need some more quantitative justification here, perhaps a reference is needed). Let the steady-state solution be $\mathbf u_{\text{true}}$. We conclude the convergence results of the velocity fields obtained from simulations on different grids in Table \ref{table:convergence_vel} and the convergence results of the boundary condition are given in Table \ref{table:convergence_bc}.
\begin{table*}[!t]
\centering
\begin{tabular}{| c| c| c|c|}
\hline
N & $||\mathbf u - \mathbf u_{\text{true}}||_1$ & $||\mathbf u - \mathbf u_{\text{true}}||_2$ & $||\mathbf u - \mathbf u_{\text{true}}||_\infty$ \\ \hline
$16$  & $2.14 \times 10^{-2}$ & $2.17 \times 10^{-2}$ & $2.22 \times 10^{-2}$\\ \hline
$32$  & $1.11 \times 10^{-2}$& $1.12 \times 10^{-2}$ & $1.13 \times 10^{-2}$\\ \hline
$64$  & $5.66 \times 10^{-3}$& $5.68 \times 10^{-3}$ & $5.68 \times 10^{-3}$\\ \hline
$128$ & $2.86 \times 10^{-3}$& $2.86 \times 10^{-3}$ & $2.85 \times 10^{-3}$\\ \hline
$256$  & $1.43 \times 10^{-3}$& $1.43 \times 10^{-3}$ & $1.08 \times 10^{-3}$\\ \hline
 \hline
\end{tabular}
\caption{Convergence results of the velocity field $\mathbf u$}\label{table:convergence_vel}
\end{table*}

\begin{table*}[t!]
\centering
\begin{tabular}{| c| c| c|c|}
\hline
$N$ & $||\mathbf{d}||_1 $& $||\mathbf{d}||_2$ & $||\mathbf{d}||_\infty$ \\ \hline
$16$  & $1.25 \times 10^{-6}$ & $1.25 \times 10^{-6}$ & $1.25 \times 10^{-6}$\\ \hline
$32$  & $6.26 \times 10^{-7}$& $6.26 \times 10^{-7}$ & $6.26 \times 10^{-7}$\\ \hline
$64$  & $3.13 \times 10^{-7}$& $3.13 \times 10^{-7}$ & $3.13 \times 10^{-7}$\\ \hline
$128$ & $1.56 \times 10^{-7}$& $1.56 \times 10^{-7}$ & $1.56 \times 10^{-7}$\\ \hline
$256$  & $7.83 \times 10^{-8}$ & $7.83 \times 10^{-8}$ & $7.83 \times 10^{-8}$\\ \hline
 \hline
\end{tabular}
\caption{Convergence results of the boundary condition}\label{table:convergence_bc}
\end{table*}
 Note that different norms ($L^1$,$L^2$,$L^\infty$) of the numerical errors are close to each other, which implies that the numerical errors are almost uniform in space. This numerical example justifies the numerical scheme described in Section 2.2 and confirms that our mesh refinement strategy deduced from the theory works.

\subsection{An application of the target point method}
The target point method described in Section 2.2 is powerful in practice and here we would like to provide a 2D example simulated by the method. 

\begin{figure}[H]
\centering
\includegraphics[width=0.9\textwidth]{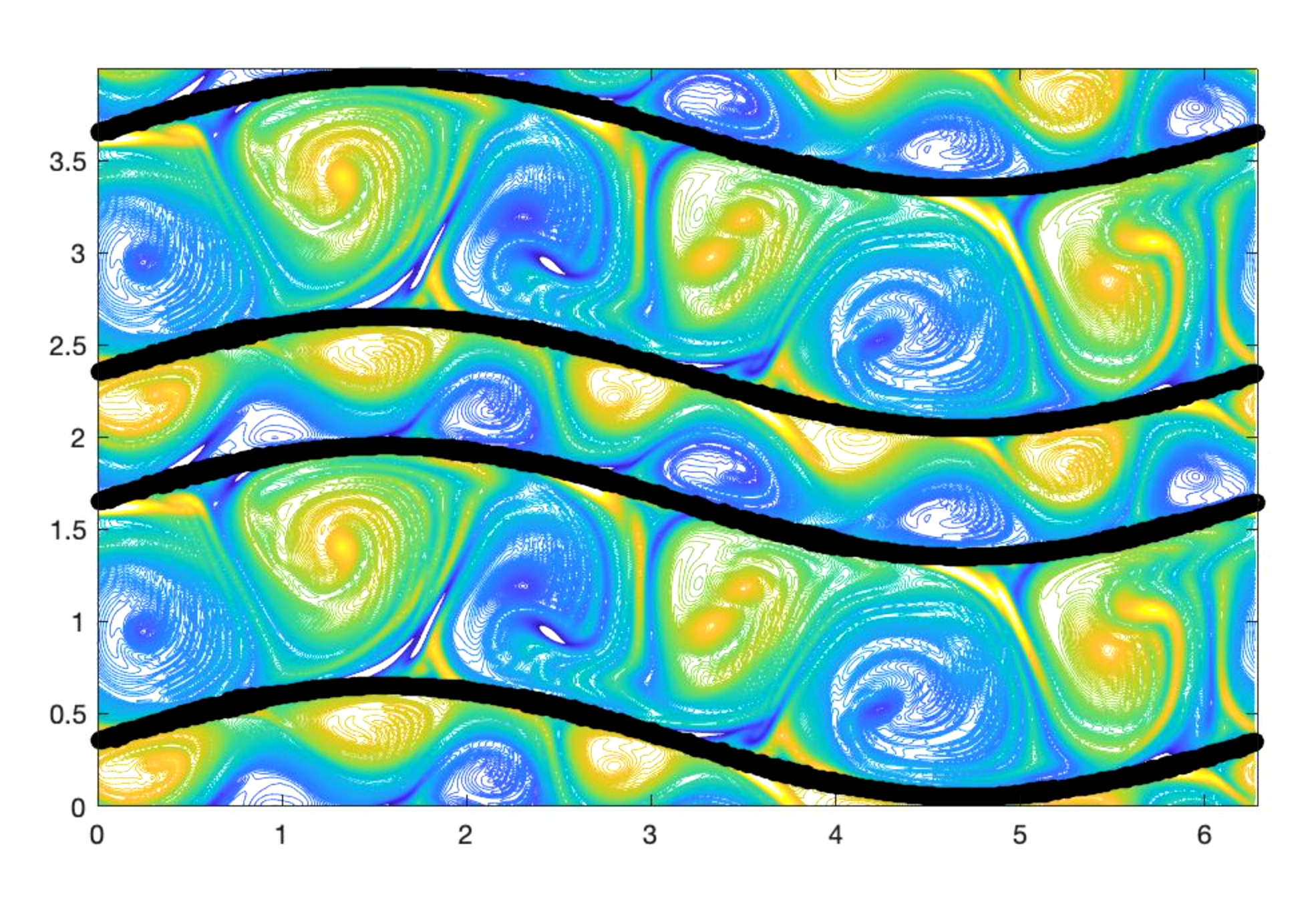}
\caption{Vorticity contour plot of stable and almost periodic flows in two channels. Two periods are shown in the $x_2$-direction for better visualization. For further details, see:
\href{https://drive.google.com/file/d/1-W5CVvs2cAVDiBPiYnwXXSEWUXA0EaK0/view?usp=sharing}{https://drive.google.com/file/d/1\-W5CVvs2cAVDiBPiYnwXXSEWUXA0EaK0/view?usp=sharing} . Please be patient while watching this video; it takes time for the interesting flow to develop, and its development as it occurs is well worth watching!}
\label{fig:6}
\end{figure}

We simulate a 2D flow in a pair of side-by-side sinusoidal
channels of different widths. The domain of our simulation is $\Omega = [0,2\pi] \times [0,2]$
 with periodic boundary conditions in both directions.
The flow starts from rest, and is
driven by a body force $\mathbf{f}^0$, which is constant in both space and time,
pointing in the $x_1$ direction. The body force is applied
everywhere, including the locations occupied by the walls. We enforce the no-slip condition on the walls by using target points using the target point method. The target points are placed at $x_2 = 0.3\sin{x_1} + 0.35$ and $x_2 = 0.3\sin{x_1} + 1.65$ for two walls respectively. We use a $512 \times 512$ Cartesian mesh for the computational domain and $4096$ target points, which are equally spaced in x, to discretize each of the two walls. Since the computational domain is not square, the mesh widths in the $x_1$ and $x_2$ directions are not equal. Whenever we refer to the \textit{width} of the channel, we mean its width in the $x_2$ direction.  Thus, the widths of our channels are $1.3$ and $2.0 - 1.3 = 0.7$ because of periodicity. All the physical and numerical parameters are listed in Table \ref{table:parameters}. 
\begin{table*}[!t]
\centering
\begin{tabular}{| c  c |c c}
$\text{Physical Parameters}$ & $\text{Values}$ &$\text{Numerical Parameters}$ & $\text{Values}$\\
$\rho$ & $1$ & $h_1$ & $2\pi/512 $\\
$\mu$ & $0.0071$ & $h_2$ & $2/512$ \\
$\mathbf{f}^0$ & $10.65$ & $\Delta t$ & $1e-5$\\
$t_{\text{total}}$ & $100$ & $K$ & $4e10$\\
$A$ & $0.3$ & $\Delta \theta$ & $2\pi/2048$\\ 
\\
$\text{Channel}$ : \\
$\text{streamwise period}$ & $2\pi$ \\
$\text{width}$ & $0.7, 1.3$ \\
\end{tabular}
\caption{The physical and numerical parameters used in simulation}\label{table:parameters}
\end{table*}
The flow accelerates under the influence of the constant body force $\mathbf{f}^0$, and at early times the vorticity field seems to be evolving
towards that of a steady flow with something like a parabolic velocity
profile, with vorticity contours being roughly parallel to the
sinusoidal channel walls.  Then, at a later time, boundary layer
separation at discrete locations becomes apparent, and this leads
rather abruptly to the formation of prominent vortices of alternate
sign that seem to fill the channel.  The vortices march downstream at
what must be the mean streamwise velocity of the fluid.  The flow then
becomes periodic in time (see below) and resembles a traveling wave,
although it cannot be strictly a traveling wave, since there is
inhomogeneity in space because of the sinusoidal channel walls. 
(Consider, for example, the curvature of the sinusoidal channel,
which is certainly not constant, so the flow is encountering
different conditions at different locations.)  A
snapshot of the voricity field is shown in Figure \ref{fig:6}.  In this figure,
and even more so in the corresponding movie, it is clear that voricity
is continually being generated at the boundaries and transferred via boundary layer
separation to
the discrete vortices that fill the channel.  The points of boundary layer separation march downstream
along with the vortices.  The vorticity that is shed from the boundary
layer is rolled up into the vortices and must ultimately be dissipated
primarily within the intense core of each vortex.

\begin{figure}[H]
\centering
\includegraphics[width=0.9\textwidth]{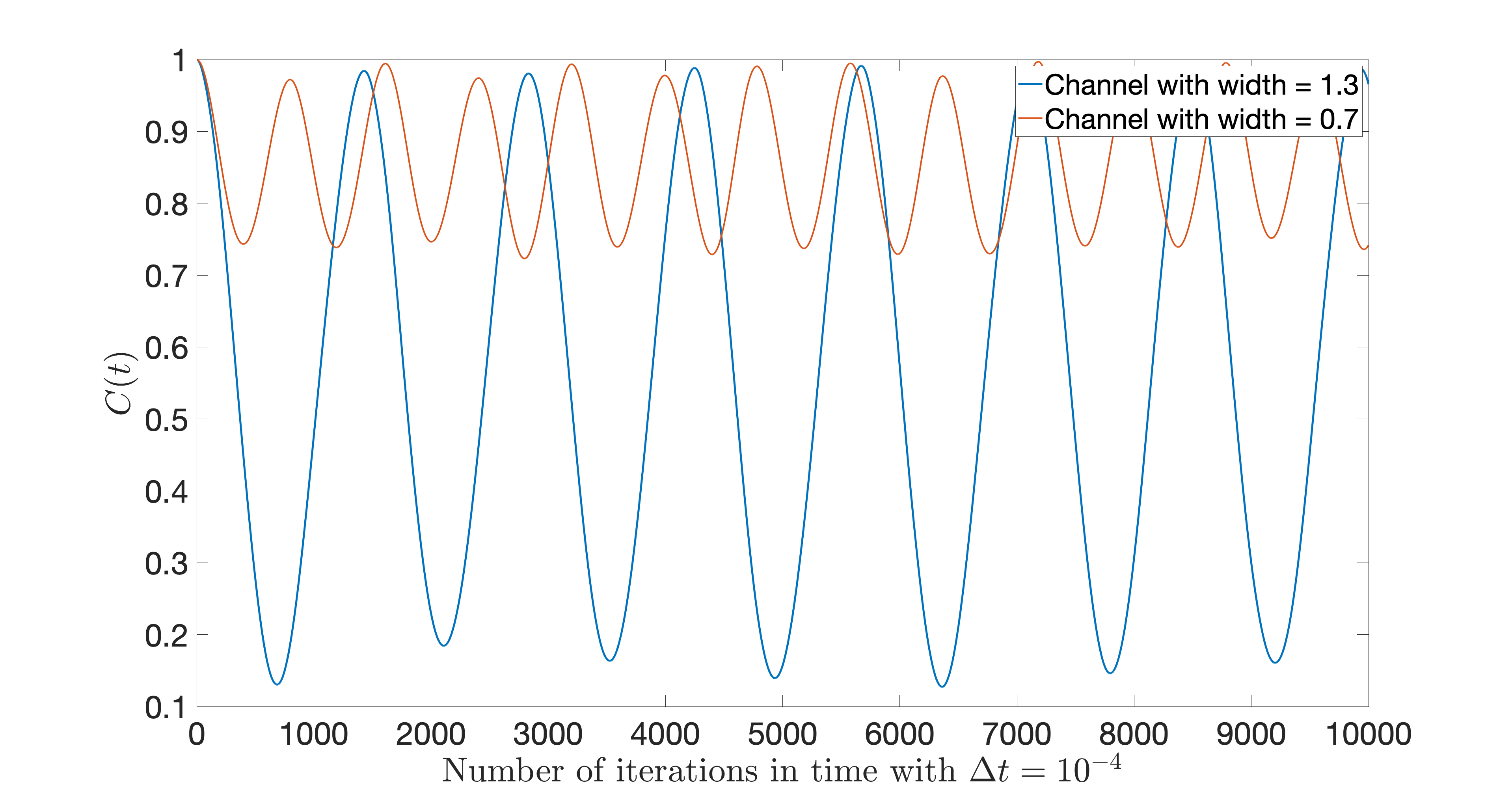}
\caption{Plot of $C(t)$ versus the number of iterations in time with $\Delta t = 10^{-4}$. The blue curve is $C(t)$ for the channel with width = 1.3 and the orange curve is $C(t)$ for the channel with width = 0.7. The flows in the two channels are both very nearly periodic with different periods. }
\label{fig:7}
\end{figure}

We use the following expression C(t) to determine
the periodicity in time of the fluid velocity field.
\begin{equation}
    C(t) = \frac{\sum_{\mathbf{x}} \mathbf{u}(\mathbf{x},t) \cdot  \mathbf{u}(\mathbf{x},t_0)}{\sqrt{\sum_{\mathbf{x}} |\mathbf{u}(\mathbf{x},t_0)|^2} \sqrt{\sum_{\mathbf{x}} |\mathbf{u}(\mathbf{x},t)|^2}}
\end{equation}
where $t_0$ is a reference time at which the flow has already become almost periodic. C(t) is the cosine of the angle (in function space) between the velocity field at time t and the velocity field at the time $t_0$. Note that we calculate $C(t)$ for each of the two channels and the sum over $\mathbf{x}$ is restricted in each case
to the grid points that are inside of the channel . If the flow is periodic for $t > t_0$, then $C(t)$ should also be periodic with peaks equal to $1$. A plot of $C(t)$ for the two channels shown in Figure \ref{fig:7}. We find that the flow is almost perfectly periodic and the wider channel has a longer period.


\subsection{Elastic membrane}
We also perform some numerical tests to show that the stability boundary predicted by the theory is in good accordance with the numerical observations for the case of a planar elastic membrane. 

\begin{figure}[H]
\centering
\includegraphics[width=0.9\textwidth]{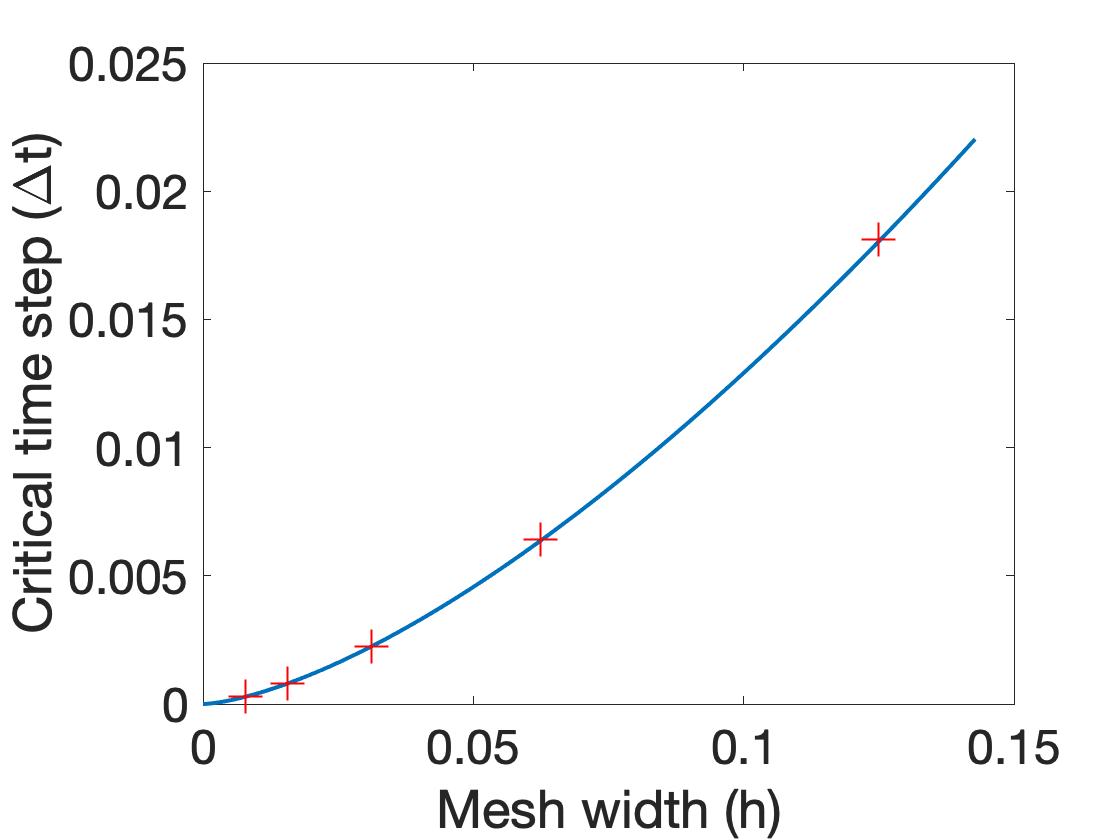}
\caption{Critical time steps for different mesh widths with $K = 100$, $\rho$ = 1, and $P = 1$. The blue curve is what we would expect based on \eqref{stability_criteria_2} and red markers are data points we get from numerical tests. The initial test velocity field is drawn i.i.d. from $0.03\mathcal{N}(0,1)$, where $\mathcal{N}(0,1)$ is the standard Gaussian distribution. For each set of parameters, we perform 10 simulations and choose the average critical time step.}
\label{fig:3}
\end{figure}

Figure \ref{fig:3} shows that the critical $\Delta t$ is proportional to $h^\frac{3}{2}$ for fixed $\rho,K,P$ when the initial velocity field is small and when we use the Navier-Stokes equations in the simulation. As we linearize the Naiver-Stokes equations in the analysis, it is not surprising that the results are good for small amplitude vibrations. Indeed, the critical $\Delta t$ only differs for at most $3\%$ when the amplitude of the Gaussian distribution we use to generate initial velocity fields increase from $0.03$ to $3$, which is almost the largest magnitude we can choose without violating the CFL condition. 

\subsection{A numerical example of an immersed elastic membrane}
\begin{figure}[h]
\centering
\includegraphics[width=0.6\textwidth]{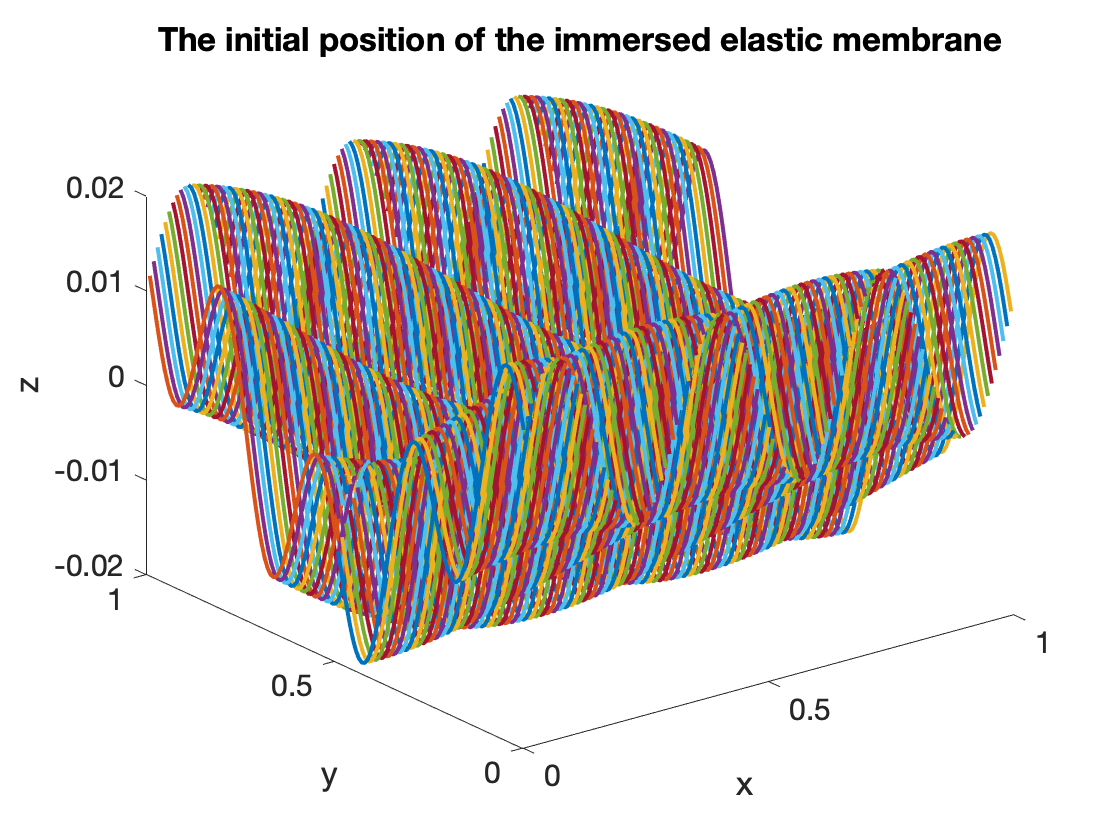}
\caption{The initial position of the immersed elastic membrane with $N = 64$ and $P = 2$. Note the exaggeration of the vertical scale in this and the following figures.}
\label{fig:membrane_initial}
\end{figure}

We simulate a 3D flow in a periodic cube with the same domain size as what we analyze in Section 2.4 and we initialize the position of the immersed elastic membrane as follows
\begin{equation}
    \mathbf{X}^0 (x_1,x_2) = (x_1.x_2,0) + A (0,0,\sin \left(2\pi (3x_1+4 x_2)\right) +\cos (2\pi x_2))
\end{equation}
where $A = 0.01$ is the amplitude of the perturbation we apply to the planar membrane. A plot of the initial position of the elastic membrane is given in Figure \ref{fig:membrane_initial}.

\begin{figure}[h]
\centering
\includegraphics[width=0.6\textwidth]{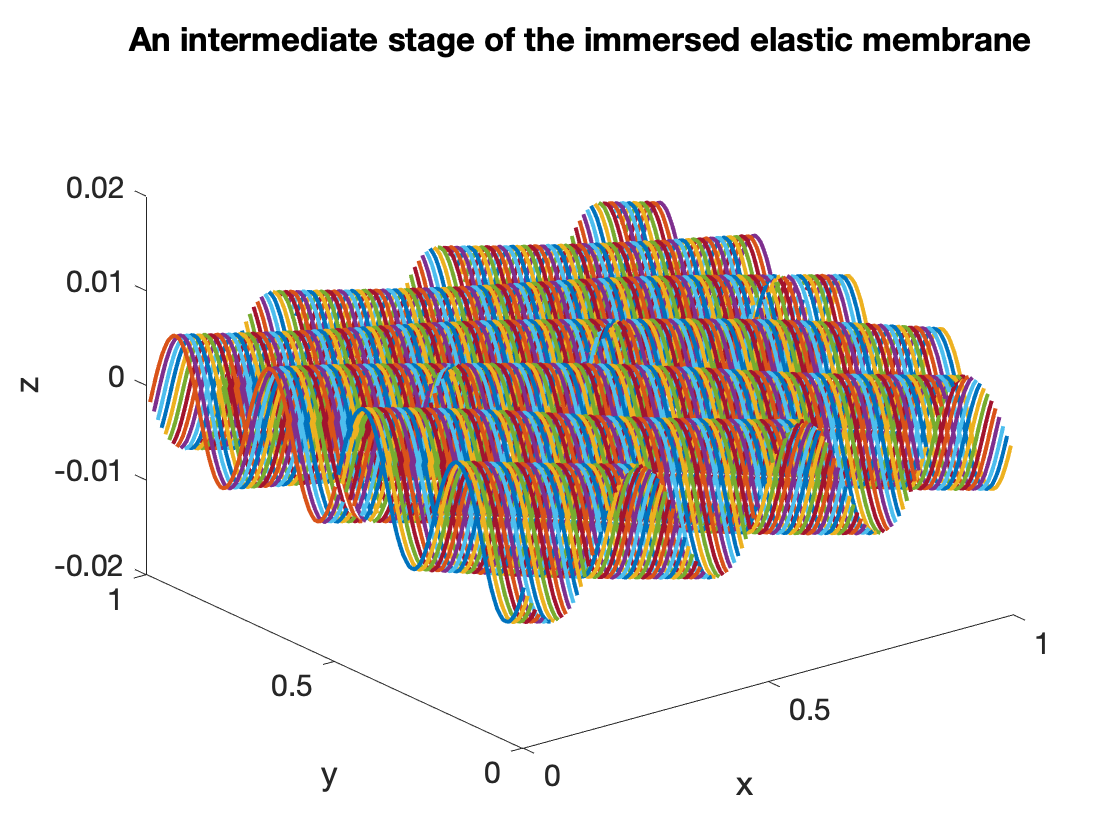}
\caption{The position of the immersed elastic membrane at time $t = 1.59 \times 10^{-2}$. For further details, see: \href{https://drive.google.com/file/d/1GX5DnCho8fDSdbQkb6wGMppMRZt8LpH-/view?usp=sharing}{https://drive.google.com/file/d/1GX5DnCho8fDSdbQkb6wGMppMRZt8LpH\-/view?usp=sharing} . In the simulation video, we shift the grid in z by $0.5$ for better visualizations and also plot the normal components of the vorticity on the three planes: $x = 0.5$, $y = 0.5$, $z = 0.5$.}
\label{fig:membrane_inter}
\end{figure}

The initial velocity field is zero and we do not apply any body force to the flow. We let the elasticity constant of this membrane be $K= 100$. After adding the small vibration to the planar elastic membrane at $z = 0$, the membrane first vibrates and then the vibration decays due to the effect of viscosity. Finally, the membrane settles down to a planar configuration at $x_3 = 0$. We also give a plot of the membrane at an intermediate time in Figure \ref{fig:membrane_inter}.

When putting the immersed boundary points, we let $P = 2$, $h = 1/64$, and the predicted critical time step from equation \eqref{stability_criteria_2} is $7.774 \times 10^{-4}$ and the actual critical time step we get for this setup is $7.771 \times 10^{-4}$. Again, we use the bisection method to determine the critical time step in the numerical experiments. 

Note that unlike what we do in the stability analysis, in the computations, we evaluate the IB delta functions at the moving positions of the immersed boundary instead of their initial positions. Moreover, we include the nonlinear terms of the Navier-Stokes equations in the computation, although they are omitted in the theory. Despite these differences, the stability boundary is well predicted by the theory. 

\begin{figure}[H]
\centering
\includegraphics[width=0.6\textwidth]{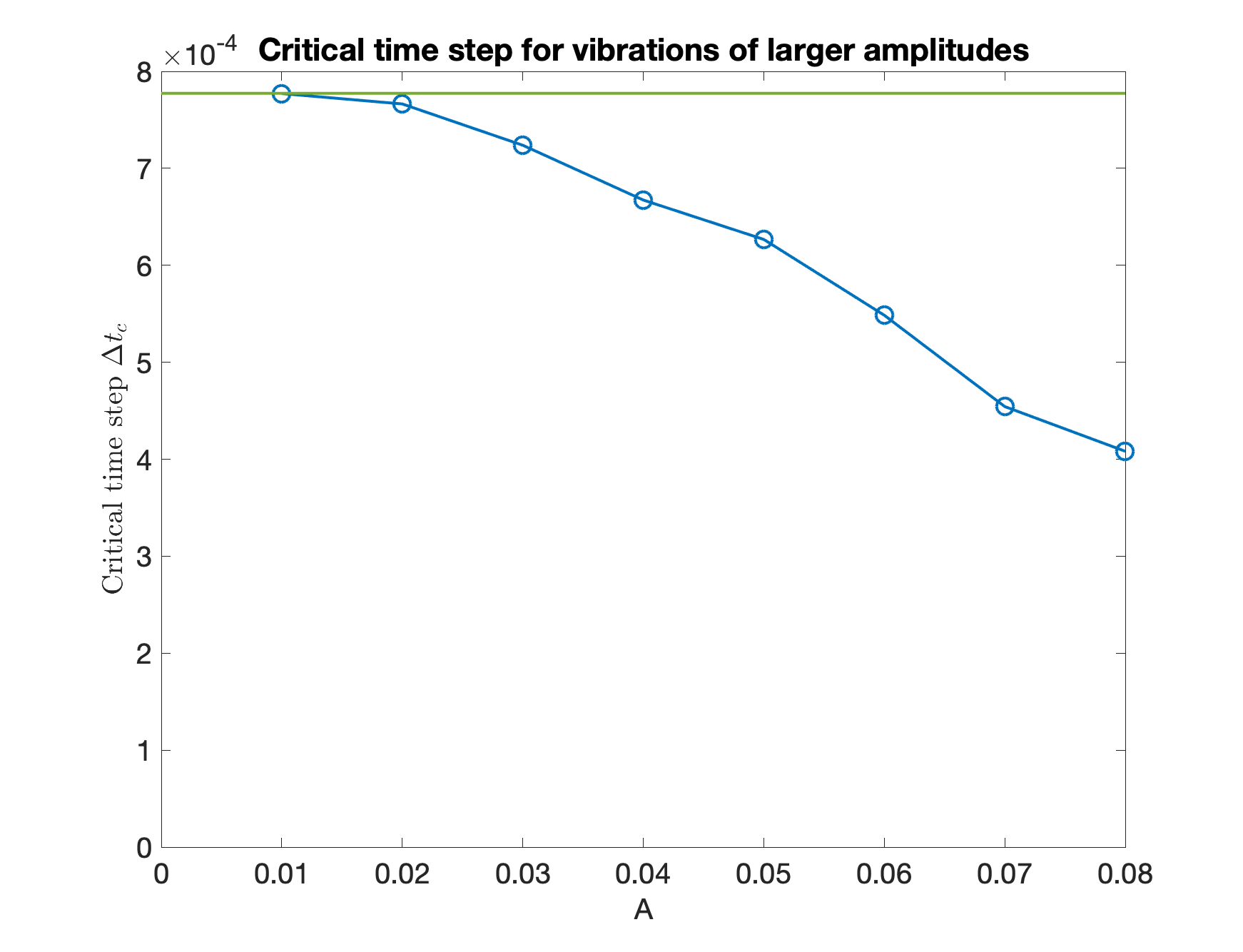}
\caption{The observed critical time step as a function of the amplitude of perturbation for $K = 100$, $h = 1/64$, $P = 2$. The limit of the blue curve as the amplitude of perturbation goes to zero is the critical time step predicted by the theory, which is indicated by the green line. }
\label{fig:membrane_relation}
\end{figure}

As the amplitude of perturbation increases, the actual critical time step deviates more from the prediction. We increase the value of $A$ and study the relation between the amplitude of perturbation and the change of critical time step. We give a plot of this relation in Figure \ref{fig:membrane_relation}.

\section{Summary and Conclusions}
In this paper, we have used Fourier analysis to study the stability of
the immersed boundary (IB) method.  To make Fourier analysis applicable,
we have considered an immersed boundary with a planar undeformed configuration,
and we have linearized the problem by considering only small-amplitude
motions of such a boundary.  The small-amplitude limit has two simplifying
effects --- one is that we can neglect the nonlinear terms in the Navier-Stokes
equations, and the other is that we can apply the boundary condition at the
undeformed location of the boundary (as in the theory of small amplitude
water waves, see for example \cite[Chapter~2]{Stoker}).  The latter simplification means
in particular that the regularized delta functions of the IB method remain
centered at fixed locations even though the boundary is in motion.  We
recommend this approach in practice as well as in theory when the goal is
to model a fixed boundary.

Our primary focus has been the use of target points to model a fixed,
no-slip boundary.  In this kind of application, the points that mark
the boundary are held in place by \textit{stiff} springs, so it is
often thought that an excessively small timestep will be required to
achieve numerical stability.  What we show in this paper is that the
timestep restriction is by no means prohibitive.  Indeed, it takes the
form $K(\Delta t)^2/h \le \text{constant}$, where $\Delta t$ is the
timestep and $h$ is the meshwidth, and it follows that we can make
$\Delta t$ be proportional to $h$ and make the stiffness parameter $K$
be proportional to $1/h$ while maintaining stability as $h\rightarrow
0$.  In this way, we simultaneously approach a continuum limit
\textit{and} it is one that obeys the desired boundary condition.  An
important remark is that the parameter $K$ in the above formula is not
the stiffness of each discrete spring.  Rather it is is the continuum
stiffness of the boundary, i.e., the force per unit area divided by
the displacement that produces that force.  Thus, $K$ has units of
force/volume.

The stability analysis that we have done is of the full IB method,
including spatial discretization, as applied to a special case.  In
particular, we allow for the boundary grid to be arbitrarily shifted
with respect to the fluid grid (although the two grids are still required
to be parallel to each other), and moreover the meshwidths of these
two grids can be different, although we do require that the boundary
meshwidth $h_B$ must be related to the fluid meshwidth $h$ by $h_B=h/P$,
where $P$ is an integer.  The analysis that we do is exact for the
case $P=1$, but for $P>1$ we obtain approximate results by making use
of the \textit{bandlimited approximation}, which we justify by
evaluating (for the first time, to our knowledge) the Fourier
transform of the standard 4-point IB delta function.  A further
justification of this approximation is that the resulting prediction
of the stability boundary is satisfied to high accuracy in numerical
experiments.  It is striking in the present work how successful the
bandlimited approximation is, and we therefore believe that it may
have future use in the analysis of the IB method.

The method of this paper is applicable to any elasticity model for the
immersed boundary itself, provided that the material of the immersed
boundary is spatially homogeneous.  To illustrate this, we have also
considered the case of an immersed membrane.  Here, the stiffness of
the membrane is a physical parameter, so it should be held constant
as the numerical parameters are refined.  The stability restriction in
the membrane case is of the form $\Delta t \le \text{constant} \;
h^{3/2}$, where the constant depends on the membrane stiffness and
also on the parameter $P$ that relates the boundary meshwidth to the
fluid meshwidth.

It is an open problem to extend the results of this paper to immersed
boundaries that are not necessarily planar, to immersed boundaries
that are undergoing large-amplitude motions, and to cases in which the
nonlinear terms of the Navier-Stokes equations play a significant role
in the dynamics.  We have, however, provided numerical evidence that
the results obtained herein are still approximately correct in such
situations.  It is our hope, therefore, that the present work will
provide a useful guide to people who use the immersed boundary method
as to what can be expected in terms of numerical stability, and at
the same time that this work can serve as inspiration for further
development of the theory of the IB method.

\section{Acknowledgement}
This work was supported in part by the National Science Foundation under grants CBET-1706562 and DMS-1646339. The authors would like to thank Sinan Gunturk for helpful discussions in relation to the the bandlimited approximation. 

\bibliographystyle{siam}
\bibliography{reference}

\end{document}